\documentclass[12pt]{article}

\usepackage{amssymb}
\usepackage[cp1251]{inputenc}   
\usepackage[english, russian]{babel}
\usepackage{amsmath}
\usepackage{amsfonts,amssymb}
\usepackage[mathcal]{eucal}

\newtheorem{theorem}{Теорема}

\newtheorem{corollary}{Следствие}
\newtheorem{proclaim}{Предложение}

\newtheorem{definition}{Определение}
\newtheorem{example}{Пример}
\newtheorem{remark}{Замечание}

\newcommand{\bo}{\hfill {$\Box$}}
\newcommand{\No}{No\ }

\renewcommand{\leq}{\leqslant}
\renewcommand{\geq}{\geqslant}
\newcommand{\cp}{ {\scriptstyle \complement}}
\newcommand{\dd}[2]{\cp^{(#2)}[{#1}]}
\newcommand{\rav}{\stackrel{\triangle}{=}}
\newcommand{\ravref}[1]{\stackrel{(\ref{#1})}{=}}
\newcommand{\leqref}[1]{\stackrel{(\ref{#1})}{\leq}}

\newcommand{\ssp}[2]{\left[#1\right]_{#2}}
\newcommand{\sspk}[2]{\Big[#1\Big]_{#2}}
\newcommand{\sspl}[1]{\Big[#1}
\newcommand{\sspll}[1]{\bigg[#1}
\newcommand{\ssplll}[1]{\Bigg[#1}
\newcommand{\sspr}[2]{#1\Big]_{#2}}
\newcommand{\ssprr}[2]{#1\bigg]_{#2}}
\newcommand{\ssprrr}[2]{#1\Bigg]_{#2}}

\begin{document}

\title{Об асимптотиках цен на больших промежутках\thanks{Институт математики и механики
		им. Н.Н.Красовского
		Уральского отделения Российской Академии Наук,
		620990, Екатеринбург,
		С.Ковалевской, 16
		Россия}}

\author{Д.В.Хлопин\\
	{\it khlopin@imm.uran.ru}}



\maketitle

\begin{abstract}
  В работе рассматриваются динамические антагонистические игры двух лиц.
    Для игр с одной и той же  динамикой, мгновенной полезности, возможностях игроков, исследуется зависимость их цен от платежной функции. Каждая платежная функция при этом понимается как мгновенная полезность, усредненная  в силу того или иного вероятностного распределения на полуоси.
   Показывается, что при выполнении принципа динамического программирования имеет место  теорема тауберова типа: из существования равномерного предела цен (при стремлении параметра масштабирования к нулю) для экпоненциальных и/или равномерных распределений следует, что тот же предел имеет место (при стремлении параметра масштабирования к нулю) с произвольной кусочно-непрерывной плотностью.   Доказан вариант такой теоремы для игр в дискретным временем. Полученные общие результаты применяются  к различным игровым постановкам, как для детерминированного, так и для стохастического случая. Также иследуются асимптотики цен при стремлении горизонта планирования к бесконечности.

\ \\

This paper is concerned with two-person dynamic zero-sum games. Let games for some family have common dynamics, running costs  and capabilities of players, and let these games differ in densities only. We show that the Dynamic Programming Principle directly leads
to the General Tauberian Theorem---that
 the existence of a uniform limit of the value functions
 for uniform distribution or for exponential distribution implies
that the value functions  uniformly converge to the same limit
 for arbitrary distribution     from large class.
 No assumptions on strategies are necessary. Applications to differential games and
 stochastic statement are considered.
 
 {\bf MSC2010}  91A25, 49L20, 49N70, 91A23, 40E05
 
 {\bf Keywords:}
 {Dynamic programming principle,  Abel mean, Cesaro mean,
 	differential games, zero-sum games}
 
 
\end{abstract}


\date{}

\maketitle
\section*{Введение}

В данной работе исследуются асимптотические свойства функции цены
для динамических игр  двух лиц при усреднении  платежа в силу заданной дисконтирующей
функции (плотности) с все меньшим  масштабом по времени.

В задачах динамической оптимизации достаточно часто, например из-за неоднозначности  момента окончания,
платеж приходится нормировать, усредняя по времени в силу того или иного вероятностного распределения.
 Чаще всего потенциальную бесконечность временного промежутка  эмулируют
  рассматривая задачи, с платежом, усредненным  по
  всё большему промежутку $[0,T]$, или усредненным с всё меньшей ставкой дисконтирования $\lambda$.
  При этом всякой реализации процесса,  функции $t\mapsto z(t)$,  в качестве платы сопоставляется интеграл
$$\int_0^\infty  \varrho(t)g(z(t))\,dt$$ от ограниченной функции мгновенной полезности $g$
в силу той или иной
дисконтирующей функции, плотности $\varrho$, после чего исследуется асимптотическое
  поведение функции цены для платежей
\begin{equation}\label{lambda}
   \int_0^\infty  \lambda\varrho(\lambda t)g(z(\lambda t))\,dt
\end{equation}
при $\lambda\to +0$.

Отметим, что чаще всего  применяют  плотности равномерного и экспоненциального распределений, получая усреднение по Чезаро, и усреднение по Абелю,
соответственно.
Теоремы, исследующие связь между  асимптотиками усреднений по Чезаро и по Абелю, с легкой руки Харди, стали называть тауберовыми  (иногда теоремами Таубера-Абеля). Не претендуя на полноту, отметим, что методы на основе таких теорем  прекрасно зарекомендовали себя в теории вероятности, теории чисел, комплексном анализе; ограничимся здесь лишь ссылками \cite{poed,ms2005,Bingham,Korevaar}. В данной работе нас будут интересовать лишь тауберовы теоремы для оптимизационных и игровых постановок.

 Существование предела цен в играх с платежами \eqref{lambda} при $\lambda\to+0$ означает,
 что цена игры малочувствительна к выбору масштабирующего параметра,
 лишь бы он был достаточно мал. В частности, именно это значение
 (asymptotic value)
 принято в стохастических постановках считать   ценой игры при бесконечном
 горизонте планирования.
 В таких постановках достаточно часто (см.~\cite{MN1981}) удается также
 найти асимптотически оптимальную стратегию,
 обеспечивающую при любом достаточно малом параметре масштабирования
 платеж, близкий к оптимальному  (uniform value), но в данной работе вопрос построения асимптотически оптимальной стратегии не рассматривается.

 Для самых разных игровых постановок удается найти условия, гарантирующие существование (при $\lambda\to +0$) равномерных пределов цен в играх,  платежи \eqref{lambda} которых усреднены в силу
 равномерного и/или экспоненциального распределений, ограничимся здесь лишь ссылками
 \cite{CQ2015,mexic,QG2013,Lions_unpub,Ziliotto2016}.
 Более того, достаточно часто
  имеет место следующая тауберова теорема: равномерная сходимость цен для платежей  \eqref{lambda}
   с равномерным и/или экспоненциальным распределением гарантирует равномерную сходимость цен к тому же пределу и для другого распределения. Такой результат был доказан для стохастических игр с конечным числом состояний и действий~\cite{MN1981}, для задач управления с дискретным временем
\cite{Lehrer},
 для задач управления общего вида~\cite{barton}, для дифференциальных игр \cite{MTIP},    для широкого класса стохастических игр \cite{Ziliotto2013}. Позже удалось  вывести \cite{arxiv2016} такую тауберову теорему во всех антагонистических играх как следствие принципа динамического программирования.

  Вариант тауберовой теоремы для  плотностей общего вида рассматривался в~\cite{1993}.
  Для детерминированных управляемых процессов с дискретным временем
 было показано, что из равномерной сходимости цен для средних по Чезаро  следует равномерная сходимость к тому же пределу и для цен с платежами  \eqref{lambda} при любой невозрастающей плотности  $\varrho$.
 Позже, для ряда стохастических постановок,
 при достаточно сильных дополнительных предположениях на динамику и возможности игроков, было показано, что в случае,
если равномерно сходятся цены для платежей,  усредненных по Чезаро, то тот же предел цен имеет место при любом выборе последовательности плотностей,  лишь бы только полная вариация плотностей стремилась к нулю \cite{Renault2013,RenaultVenel2012,Ziliotto2016}. Для детерминированных постановок в эргодическом случае  не менее сильные результаты показаны в \cite{CQ2015,LiQixRenault}.

 Таким образом, достаточно часто удается показать,
  что из равномерной сходимости цен для  платежей,  усредненных по Чезаро, следует асимптотическая нечувствительность цен к выбору для платежа \eqref{lambda} плотности $\varrho$ из достаточно широкого класса.
 Основная цель данной работы --- вывести такую тауберову теорему напрямую из принципа динамического программирования, не требуя каких-либо технических предположений на динамику, возможности игроков, существования в этих играх седловой точки и т. п. Для равномерного и экспоненциального распределений такая тауберова теорема была показана в \cite{arxiv2016}.

 Сама работа построена следующим образом. Сначала вводится абстрактная игровая постановка: следуя \cite{barton,arxiv2016} постулируются множества позиций и процессов, функция мгновенной полезности, затем, как и в \cite{arxiv2016}, собственно игра сводится к одному отображению, ``{\it игровому отображению}'', переводящему каждый ограниченный функционал над множеством процессов, ``{\it платежную функцию}'', в некоторый ограниченный функционал над множеством позиций, ``{\it функцию цены}''. Там же вводятся плотности, операции над плотностями
   $\varrho^{T}_{\textrm{shift}}$, $\varrho^\lambda_{\textrm{scale}}$ и наконец постулируется принцип динамического программирования. В параграфе~\ref{main}  мы формулируем общий результат, тауберову теорему для игрового отображения, как в непрерывном (теорема~\ref{4}), так и в дискретном (теорема~\ref{44}) времени.  В следующем параграфе, подставляя те или иные игровые отображения, мы получаем соответствующие тауберовы теоремы для игр в нормальной форме, для дифференциальных игр, для стохастических постановок. Подробное обсуждение предположений  теорем~\ref{4} и~\ref{44}
   проводится в параграфе~\ref{prooff}, сразу после  ввода ряда вспомогательных определений и утверждений.
   Оставшаяся часть работы, параграфы \ref{proofg}--\ref{doctochi}, посвящена доказательству этих
   теорем.

\section{Общая постановка}\label{general}

 Знак ``$\rav$'' всюду далее будет означать ``равенство по
 определению''. Примем также $\mathbb{R}_+\rav\mathbb{R}_{\geq 0}.$

\subsection*{Динамическая система}
 Всюду далее будем считать заданными:
 \begin{itemize}
   \item   некоторое непустое множество  $\Omega$ (множество состояний);
   \item   некоторое непустое множество $\mathbb{K}$ отображений из $\mathbb{R}_+$ в $\Omega$ (множество процессов);
   \item   некоторое отображение $g:\Omega\mapsto [0,1]$\ (функция мгновенной полезности).
 \end{itemize}
 Будем также предполагать, что   для всех процессов $z\in \mathbb{K}$ суперпозиции $\mathbb{R}_+ \ni t\mapsto g(z(t))$ измеримы по Лебегу.

\subsection*{ О платежах, ценах и игровом отображении}

  Обозначим через  $\mathfrak{C}$ и  $\mathfrak{U}$
  множества всех ограниченных отображений в $\mathbb{R}$ из, соответственно, $\mathbb{K}$ и $\Omega$.
 Фактически множество   $\mathfrak{C}$ используется как множество
всех возможных (для данной динамики) платежей, а множество
  $\mathfrak{U}$ заведомо содержит все цены для игр с платежами из $\mathfrak{C}.$

Отображение $\mathbb{V}$ из $\mathfrak{C}$ в $\mathfrak{U}$ назовем {\it игровым отображением}, если имеет место:
\begin{subequations}
\begin{eqnarray}
\label{conditions1}
  \mathbb{V}[Ac+B]=A\,\mathbb{V}[c]+B\qquad  \forall  c\in\mathfrak{C}, A\geq 0, B\in\mathbb{R},\\
  \Big(c_1(z)\leq c_2(z)\ \forall
  z\in\mathbb{K}\Big)\Rightarrow\Big(\mathbb{V}[c_1](\omega)\leq \mathbb{V}[c_2](\omega)\ \forall \omega\in\Omega\Big).
\label{conditions2}
\end{eqnarray}
\end{subequations}
 Указанные соотношения более чем типичны для отображений, сопоставляющих каждому платежу соответствующую цену игры.
  Примеры таких отображений, для различных игровых постановок, смотрите далее в параграфе \ref{appli}.

\subsection*{Плотности. Операции с плотностями}
Обозначим через~$\mathfrak{D}$ множество всех плотностей $\varrho$ с
$\textrm{supp}\,\varrho\subset\mathbb{R}_+$. Далее можно считать все элементы в
~$\mathfrak{D}$ отображениями из $\mathbb{R}_+$ в $\mathbb{R}_+$.

Для любых плотности  $\varrho\in\mathfrak{D}$ и положительного числа  $T>0$, в случае, если выполнено
$\int_T^{\infty}\varrho(t)\,dt>0$, зададим плотность $\varrho^{T}_{\textrm{shift}}\in\mathfrak{D}$ правилом:
$$\varrho^{T}_{\textrm{shift}}(t)=\frac{\varrho(t+T)}{\int_T^{\infty}\varrho(t)\,dt}\qquad\forall t\geq 0.$$

Для любых плотности  $\varrho\in\mathfrak{D}$ и положительного числа  $\lambda>0$ введем
плотность $\varrho^\lambda_{\textrm{scale}}\in\mathfrak{D}$ следующим образом: $$\varrho^{\lambda}_{\textrm{scale}}(t)=\lambda\varrho(\lambda t)\qquad\forall t\geq 0.$$

\subsection*{Усредненные в силу плотности платежи}

Для любой плотности $\varrho\in\mathfrak{D}$, усредним функцию мгновенной полезности  вдоль каждого процесса. Мы получим некоторый платеж $\cp[\varrho]\in\mathfrak{C}$, как  функцию из  $\mathbb{K}$ в $\mathbb{R}$, действующую по правилу:
\begin{eqnarray*}
    \cp[\varrho](z)\rav\int_{0}^\infty \varrho(t) g(z(t))\,dt\in[0,1]\quad\forall z\in\mathbb{K}.
\end{eqnarray*}
  Теперь, фиксируя некоторое игровое отображение $\mathbb{V}$, мы  автоматически
  находим цену (в силу этого отображения) для платежа $\cp[\varrho]$, то есть цену, соответствующую  плотности $\varrho$; такую суперпозицию $\mathbb{V}\circ\cp$, действующую из $\mathfrak{D
  }$ в $\mathfrak{U}$, обозначим через $\mathcal{V}$:
\begin{eqnarray}\label{superp}
    \mathcal{V}[\varrho](\omega)\rav\mathbb{V}\big[\cp[\varrho]\big](\omega)\qquad\forall z\in\mathbb{K}.
\end{eqnarray}

Отметим две полезные в дальнейшем  оценки.
Легко видеть,
условия
\eqref{conditions1} и~\eqref{conditions2} гарантируют
\begin{equation}\label{o0}
\sup_{\omega\in\Omega}|\mathbb{V}[c_1](\omega)-\mathbb{V}[c_2](\omega)|\leq \sup_{z\in\mathbb{K}}|c_1(z)-c_2(z)|
\end{equation}
для любых платежей $c_1,c_2\in\mathfrak{C}$;
откуда для любых плотностей $\nu,\nu'\in\mathfrak{D}$ имеет место
\begin{equation}\label{o1}
\sup_{\omega\in\Omega}|\mathcal{V}[\nu](\omega)-\mathcal{V}[\nu'](\omega)|\leq \int_{0}^\infty |\nu(t)-\nu'(t)|\,dt.
\end{equation}

\subsection*{Принцип динамического программирования}
Рассмотрим некоторые плотность~$\varrho\in\mathfrak{D}$ и
положительное число $T>0$, в случае если
$\int_{T}^{\infty}\varrho(t)\,dt>0$ плотность
$\varrho^{T}_{\textrm{shift}}$ определена, и мы можем ввести платеж
$\dd{\varrho}{T}:\mathbb{K}\to\mathbb{R}$ следующим правилом: для
всех $z\in\mathbb{K}$,
   $$\dd{\varrho}{T}(z)\rav
   \int_{0}^{T}\varrho(t)g\big(z(t)\big)\,dt+\int_{T}^{\infty}\varrho(t)\,dt\cdot
   \mathbb{V}\big[\cp\big[\varrho^{T}_{\textrm{shift}}\big]\big]\big(z(T)\big).$$

\begin{definition}
Мы скажем, что игровое отображение~$\mathbb{V}$  удовлетворяет на~$\mathfrak{D}$  принципу динамического программирования, если для всех $\varrho\in\mathfrak{D},T>0$ из
 $\int_{T}^{\infty}\varrho(t)\,dt>0$ следует, что
   $\mathbb{V}\big[\dd{\varrho}{T}\big]\equiv\mathbb{V}\big[\cp[\varrho]\big],$ то есть цены игр с платежами $\dd{\varrho}{T}$ и $\cp[\varrho]$ совпадают.
\end{definition}
 Отметим, что этот принцип достаточно часто имеет место, в частности он выполнен
 в задачах управления и  в классических постановках дифференциальных игр. Тем не менее,
 для задач с дискретным временем лучше это определение ослабить.
 \begin{definition}
    Мы скажем, что игровое отображение~$\mathbb{V}$  удовлетворяет на~$\mathfrak{D}$ дискретному варианту принципа динамического программирования, если для всех $\varrho\in\mathfrak{D},n\in\mathbb{N}$ из
    $\int_{n}^{\infty}\varrho(t)\,dt>0$ следует, что
    $\mathbb{V}\big[\dd{\varrho}{n}\big]\equiv\mathbb{V}\big[\cp[\varrho]\big],$ то есть цены игр с платежами $\dd{\varrho}{n}$ и $\cp[\varrho]$ совпадают.
 \end{definition}

\section{Основной результат}
\label{main}
Введем следующие обозначения для плотностей равномерных и экспоненциальных распределений:
$$\varpi_T(t)=\frac{1}{T} \cdot 1_{[0,T]}(t),\quad \pi_\lambda(t)=\lambda\cdot e^{-\lambda t} \qquad\forall     \lambda,T>0,t\geq 0.$$
Будем говорить, что плотность  $\varrho\in\mathfrak{D}$ степенная, если для некоторых положительных чисел $\alpha,\beta,\gamma$  выполнено
$$\varrho(t)=\frac{1}{(\alpha+\beta t)^{\gamma}}\qquad \forall t\geq 0.$$


В параграфах \ref{proofg}--\ref{doctochi} будут доказаны следующие теоремы:
\begin{theorem}
\label{4}
 Предположим, что  некоторое игровое отображение $\mathbb{V}$  удовлетворяет
 на~$\mathfrak{D}$ принципу динамического программирования,
 отображение $\mathcal{V}\colon\mathfrak{D}\to\mathfrak{U}$
 задано правилом \eqref{superp}, а кроме того, некоторая функция $U_*\colon\Omega\to [0,1]$ выбрана произвольным образом.

Тогда  следующие условия эквивалентны:
 \begin{description}
 \begin{subequations}
\item[$(c)$]\
  для каждой кусочно-непрерывной на  $(0,\infty)$ плотности  $\mu\in \mathfrak{D}$
   цены $\mathcal{V}[\mu^\lambda_{\textrm{scale}}]$ сходятся к~$U_*$ равномерно
   по $\omega\in\Omega$ при $\lambda\to  +0$, то есть
    \begin{eqnarray}\label{1113b}
       \lim_{\lambda\to +0} \sup_{\omega\in\Omega}
       \Big|\mathcal{V}[\mu^\lambda_{\textrm{scale}}](\omega)-U_*(\omega)\Big|=0;
    \end{eqnarray}
   \item[$(u)$] \
   цены $\mathcal{V}[\varpi_T]$  сходятся к~$U_*$ равномерно по $\omega\in\Omega$ при $T\to+\infty$, то есть
    \begin{eqnarray}\label{1113u}
       \lim_{T\to+\infty} \sup_{\omega\in\Omega}
       \Big|\mathcal{V}[\varpi_T](\omega)-U_*(\omega)\Big|=0;
     \end{eqnarray}
  \item[$(e)$]\
   цены $\mathcal{V}[\pi_\lambda]$ сходятся к~$U_*$ равномерно по $\omega\in\Omega$ при $\lambda\to +0$, то есть
     \begin{eqnarray}\label{1113e}
       \lim_{\lambda\to+0} \sup_{\omega\in\Omega}
        \Big|\mathcal{V}[\pi_\lambda](\omega)-U_*(\omega)\Big|=0;
       \end{eqnarray}
  \item[$(p)$] \  для некоторой степенной  плотности $\varrho\in \mathfrak{D}$
       цены $\mathcal{V}[\varrho^{T}_{\textrm{shift}}]$ сходятся к~$U_*$ равномерно по $\omega\in\Omega$ при $T\to+\infty$, то есть
\begin{equation}\label{1113m}
    \lim_{T\to+\infty} \sup_{\omega\in\Omega}
  \Big|\mathcal{V}[\varrho^{T}_{\textrm{shift}}](\omega)-U_*(\omega)\Big|=0;
    \end{equation}
  \item[$(q)$] \    для всякой степенной плотности $\varrho\in \mathfrak{D}$
  выполнено \eqref{1113m}, то есть
       цены $\mathcal{V}[\varrho^{T}_{\textrm{shift}}]$ сходятся к~$U_*$ равномерно по $\omega\in\Omega$ при $T\to+\infty$.
\end{subequations}
         \end{description}

   В частности, если для некоторых
   кусочно-непрерывной на  $(0,\infty)$ плотности  $\mu\in \mathfrak{D}$ и  степенной плотности $\varrho\in \mathfrak{D}$ хотя бы один среди первых трех
   пределов в
   $$\lim_{T\to+\infty}\mathcal{V}[\varpi_{T}],\
   \lim_{\lambda\to+0}\mathcal{V}[\pi_{\lambda}],\  \lim_{T\to+\infty}\mathcal{V}[\varrho^{T}_{\textrm{shift}}],\
   \lim_{\lambda\to+0}\mathcal{V}[\mu^{\lambda}_{\textrm{scale}}]
   $$
   существует и равномерен на $\Omega,$  то все эти четыре предела существуют,
   равномерны на $\Omega$, совпадают  и не зависят от выбора кусочно-непрерывной
    на  $(0,\infty)$ плотности  $\mu\in \mathfrak{D}$ и  степенной плотности
     $\varrho\in \mathfrak{D}$.
\end{theorem}
\begin{theorem}
    \label{44}
     Предположим, что  некоторое игровое отображение $\mathbb{V}$  удовлетворяет на~$\mathfrak{D}$ дискретному варианту принципа динамического программирования, отображение $\mathcal{V}\colon\mathfrak{D}\to\mathfrak{U}$ задано правилом \eqref{superp}, а некоторая функция $U_*\colon\Omega\to[0,1]$ выбрана произвольным образом.
     Пусть также
     выполнено
    \begin{equation}
    \label{579}
    z(n+s)=z(n)\qquad\forall z\in\mathbb{K}, n\in\mathbb{N}\cup\{0\}, s\in(0,1).
    \end{equation}

    Тогда каждое из условий $(c)$, $(u)$, $(e)$, $(p)$, $(q)$ эквивалентно условию
    \begin{description}
        \item[$(u')$] \
        цены $\mathcal{V}[\varpi_{n}]$ сходятся к~$U_*$ равномерно на $\Omega$ при $n\to\infty$.
    \end{description}
\end{theorem}
Условия этих теорем будут обсуждаться в параграфе~\ref{prooff}. Сейчас мы применим эти теоремы к конкретным игровым постановкам.

\section{Тауберовы теоремы для различных игровых постановок}
\label{appli}
 \subsection{Игры в нормальной форме}

Пусть, как и раньше, заданы множества $\Omega$ и $\mathbb{K}$, а также отображение $g:\Omega\to[0,1]$.

  Для каждого $\omega\in\Omega$ рассмотрим некоторые непустые множества $\mathcal{L}(\omega)$ и $\mathcal{M}(\omega)$ правил первого и второго игроков соответственно.
Пусть для любого    $\omega\in\Omega$ каждой паре $(l,m)\in\mathcal{L}(\omega)\times\mathcal{M}(\omega)$ правил игроков
соответствует единственный процесс $\zeta[\omega,l,m]\in \mathbb{K}$.

Каждая функция платы $c\in\mathfrak{C}$ задает следующую игру:
  для всякого $\omega\in\Omega$
  первый игрок показывает некоторое $l\in \mathcal{L}(\omega)$, после чего второй игрок выбирает  $m\in \mathcal{M}(\omega)$,
  в этом случае реализуется процесс $\zeta[\omega,l,m]\in \mathbb{K}$, в ходе которого второй  платит первому  платеж $c(\zeta[\omega,l,m]).$
  Функция цены в этой игре, нижняя цена, имеет вид:
  \begin{eqnarray}
 V[c](\omega)\rav\sup_{l\in \mathcal{L}(\omega)}\inf_{m\in \mathcal{M}(\omega)}c(\zeta[\omega,l,m])\qquad \forall\omega\in\Omega.
 \label{295}
\end{eqnarray}
  Например для каждой плотности $\varrho\in\mathfrak{D}$ мы получаем:
  \begin{eqnarray*}
 \mathcal{V}[\varrho](\omega)&=&\sup_{l\in \mathcal{L}(\omega)}\inf_{m\in \mathcal{M}(\omega)}\int_{0}^\infty
 \varrho(t)
 g(\zeta[\omega,l,m](t))\,dt\qquad \forall\omega\in\Omega.
\end{eqnarray*}

\begin{theorem}\label{normal1}
Предположим, что для всех  $T>0$ и $\varrho\in\mathfrak{D}$
 значение заданного правилом \eqref{295} отображения
$V$ для платежа
  \begin{eqnarray}
 \mathbb{K}\ni z&\mapsto& \int_{0}^T \varrho(t)g(z(t))\,dt\nonumber\\
 &\ &+
 \sup_{l\in \mathcal{L}(z(T))}\inf_{m\in \mathcal{M}(z(T))}
 \int_{0}^\infty \varrho(t+T)g(\zeta[z(T),l,m](t))\,dt
 \label{407}
\end{eqnarray}
 совпадает с $\mathcal{V}[\varrho]$.
  Пусть также $U_*$ --- некоторая функция из $\Omega$ в $[0,1]$.

 Тогда результат теоремы~\ref{4} выполнен, в частности условия
$(c)$, $(u)$, $(e)$, $(p)$, $(q)$ эквивалентны.
\end{theorem}
\underline{Доказательство теоремы~\ref{normal1}.}
Легко видеть, что введенное правилом \eqref{295} отображение $V$
удовлетворяет условиям \eqref{conditions1},\eqref{conditions2}, то
есть является игровым отображением. Поскольку платеж   \eqref{407}
совпадает с $\dd{\varrho}{T}$, имеет место и принцип динамического
программирования  на
 $\mathfrak{D}$. Теперь
теорема~\ref{normal1} следует напрямую из теоремы~\ref{4}.\bo


 Для игр в нормальной форме  эквивалентность условий $(e),(u)$ в сходных предположениях  была показана в \cite[Theorem~1]{arxiv2016}. Условия на $\mathbb{K}$, гарантирующие равное $\mathbb{V}\big[\cp[\varrho]\big]$  значение отображения
 $V$ для платежа
 \eqref{407}, смотрите в
 \cite[(2.2)-(2.4)]{MatSb}, \cite{SIAM}.

В случае одного игрока (например потребовав
$\mathcal{M}(\omega)=\{\mathbb{K}\}$,
$\mathcal{L}(\omega)=\{{z\in\mathbb{K}\,|\,z(0)=\omega}\}$ для всех
$\omega\in\Omega$, как в \cite[параграф 7]{MatSb}), мы  имеем
\begin{equation}
\label{425}
  \mathcal{V}[\varrho](\omega)\rav\sup_{z\in\mathbb{K},z(0)=\omega} \int_{0}^\infty
 \varrho(t)
 g(z(t))\,dt
 \qquad\forall z\in\mathbb{K},\varrho\in\mathfrak{D}.
 \end{equation}
Применяя теперь  теорему \ref{normal1},мы  получаем следующий результат для задач управления общего вида:
   \begin{corollary}
     \label{u1}
Пусть для всех
 $\omega\in\Omega$, $\varrho\in\mathfrak{D}$ имеет место
\begin{eqnarray}
  \mathcal{V}[\varrho](\omega)=&{\displaystyle\sup_{z\in\mathbb{K},z(0)=\omega}}&
\Big[\int_{0}^T
 \varrho(t)
 g(z(t))\,dt\nonumber\\ &\ &+
\sup_{z_1\in\mathbb{K},z_1(0)=z(T)}
\int_{0}^\infty
 \varrho(t+T)
 g(z_1(t))\,dt
 \Big].\label{435}
\end{eqnarray}

 Тогда,  для заданного правилом \eqref{425} отображения $\mathcal{V}$, при любом выборе функции $U_*:\Omega\to\mathbb{R}$, имеет меcто результат теоремы \ref{4}, в частности
условия $(c)$, $(u)$, $(e)$, $(p)$, $(q)$ эквивалентны.
    \end{corollary}
  Заметим, что из следствия~\ref{u1} можно легко извлечь тауберовы теоремы для любых  задач управления детерминированными процессами.
  Подобная, столь общая формулировка  (для эквивалентности $(e)\Leftrightarrow (u)$) впервые появилась в работе  в \cite{barton}, там же имеется подробный исторический экскурс о применимости тех или иных тауберовых теорем   в играх и задачах управления с непрерывным временем. Тем не менее,
   результат следствия~\ref{u1} для конкретных классов задач управления может быть существенно усилен хотя бы в части более общих условий на семейство плотностей.
   Например, для автономных управляемых систем  в работе \cite{LiQixRenault} исследуются асимптотики цен $\mathcal{V}[\varrho_\lambda]$ лишь при следующем  условии на плотность:
    \begin{eqnarray}\label{501}
    \sup_{S>0}\lim_{\lambda\to 0} \sup_{T>0,s<S}|\varrho_\lambda(T)-\varrho_\lambda(T+s)| =0.
    \end{eqnarray}

  Условия существования равномерного предела цен в случае одного игрока исследовались во многих работах.
   Так, в работах \cite{Gaitsgori1986,ArisawaLions,AG2000,Barles2000,wKAM1998,Lions_unpub} были найдены предположения (такие как существование эргодической меры), при которых цены сходятся к константе. Пределы, отличные от констант, рассматривались в работах \cite{CQ2015,QG2013,grune,QuinRenault}.

\subsection{Дифференциальные игры}

 В некотором конечномерном евклидовом пространстве~$\mathbb{X}$ рассмотрим автономную конфликтно-управляемую систему
\begin{equation}\label{sys}
 \dot{x}=f(x,a,b),\ x(0)\in \mathbb{X},\   a(t)\in \mathbb{A},\ b(t)\in \mathbb{B};
\end{equation}
здесь,~$\mathbb{A}$ и~$\mathbb{B}$ --- некоторые непустые компактные подмножества конечномерных евклидовых пространств.

 Под множеством  программных управлений первого (соответственно второго) игрока будем понимать множество всех измеримых по
 Борелю, ограниченных функций из $\mathbb{R}_+$ в $\mathbb{A}$ (соответственно
 в $\mathbb{B}$). Обозначим множества всех программных управлений первого и второго игрока через $B(\mathbb{R}_+,\mathbb{A})$ и $B(\mathbb{R}_+,\mathbb{B}).$

 Предположим, что функции
 $f:\mathbb{X} \times \mathbb{A} \times \mathbb{B}\to\mathbb{X}$ и $g:\mathbb{X} \times \mathbb{A} \times \mathbb{B}\to [0,1]$ непрерывны,
  и для некоторой положительной константы  $L>0$ при любых $x,y\in\mathbb{X}$, $a\in \mathbb{A}$, $b\in \mathbb{B}$ выполнено
  $$\big|\big|f(x,a,b)-f(y,a,b)\big|\big|+\big|g(x,a,b)-g(y,a,b)\big|\leq L\big|\big|x-y\big|\big|.
  $$
 Теперь для всякого начального условия $x(0)=x_*\in\mathbb{X}$, каждая пара
 управлений  $(a,b)\in B(\mathbb{R}_+,\mathbb{A})\times B(\mathbb{R}_+,\mathbb{B})$ игроков задает
 единственное (определенное на всем  $\mathbb{R}_+$) решение  $x(\cdot)=y(\cdot;x_*,a,b)$ системы
 \eqref{sys}.

 В отсутствии условия седловой точки в маленькой игре (условия Айзекса)
 нижняя цена и верхняя цена игры вообще говоря зависят от выбранной для игроков
 формализации их стратегии
  \cite[глава~XVI]{ks}, \cite[Subsect.~14]{subb}. В данной работе мы применим формализацию на основе квазистратегий (cтратегий Элиотта-Калтона) \cite{EK}.

\begin{definition}
      Отображение  $\alpha : B(\mathbb{R}_+,\mathbb{B})\mapsto B(\mathbb{R}_+,\mathbb{A})$ назовем   {\it квазистратегией  первого игрока},
      если для всех   $t>0$ и $b,b'\in B(\mathbb{R}_+,\mathbb{B})$, из
      $b|_{[0,t]}   =  b'|_{[0,t]}$  следует $\alpha[b]|_{[0,t]}  =  \alpha[b']|_{[0,t]}$.
  Через $\mathcal{A}$ обозначим множество всех квазистратегий первого игрока.
\end{definition}

  Каждой плотности  $\varrho\in\mathfrak{D}$ определим для всех  $x_*\in\Omega,$
     \begin{eqnarray}
   {\mathcal{V}}[\varrho](x_*)\rav\inf_{\alpha\in\mathcal{A}}\sup_{b\in B(\mathbb{R}_+,\mathbb{B})}
   \int_0^\infty \varrho(t)g(y(t;x_*,\alpha(b),b),\alpha(b)(t),b(t))\,dt.
   \label{1588}
\end{eqnarray}

\begin{theorem}\label{diff}
      Пусть некоторое непустое множество  $\Omega\subset\mathbb{X}$ сильно инвариантно относительно системы \eqref{sys}, то есть $y(t;x_*,a,b)$ лежит в $\Omega$
            для всех $x_*\in\Omega,$ $t\geq 0$, $(a,b)\in B(\mathbb{R}_+,\mathbb{A})\times B(\mathbb{R}_+,\mathbb{B})$.

Тогда, если отображение~ ${\mathcal{V}}$ задано правилом \eqref{1588}, то при любом выборе функции $U_*:\Omega\to [0,1]$ имеет место результат теоремы~\ref{4}, в частности
условия  $(c)$, $(u)$, $(e)$, $(p)$, $(q)$ эквивалентны.
\end{theorem}
\underline{Доказательство теоремы~\ref{diff}.}
 Примем
    $\Omega'\rav\Omega\times\mathbb{A}\times\mathbb{B},\ \mathbb{K}'\rav C(\mathbb{R}_+,\Omega)\times B(\mathbb{R}_+,\mathbb{A})\times B(\mathbb{R}_+,\mathbb{B}).$
     Пусть ${\mathfrak{C}}'$ и $\mathfrak{U}'$
     множества всех ограниченных  отображений в $\mathbb{R}$  из, соответственно,  $\mathbb{K}'$ и  $\Omega'$.

 Всякому платежу  $c\in{\mathfrak{C}'}$ сопоставим функцию цены правилом:
 для всех $\omega=(x_*,a_*,b_*)\in\Omega'$
     \begin{eqnarray}\label{1587}
   {V}[c](x_*,a_*,b_*)\rav\inf_{\alpha\in\mathcal{A}}\sup_{b\in B(\mathbb{R}_+,\mathbb{B})}
   c(y(\cdot;x_*,\alpha(b),b),\alpha(b),b).
\end{eqnarray}
Легко видеть, что условия \eqref{conditions1},\eqref{conditions2} выполнены, и
 ${V}$ является игровым отображением.
 Поскольку ${V}[c](x_*,a_*,b_*)={V}[c](x_*,a'_*,b'_*)$ имеет место для всех  $(x_*,a_*,b_*)$, $(x_*,a'_*,b'_*)\in\Omega'$,
  найдется некоторое отображение $\bar{V}[c]:\Omega\to\mathbb{R}$, для которого
${V}[c](x_*,a_*,b_*)=\bar{V}[c](x_*)$.

  Осталось проверить, что  $\bar{V}$ удовлетворяет на  $\mathfrak{D}$ принципу динамического программирования.
   Для этого достаточно убедиться, что для всех $\varrho\in\mathfrak{D}$, $T>0$ из
    $\int_T^\infty\varrho(t)\,dt<1$ следует, что  $\bar{V}[{\cp}[\varrho]]\equiv \bar{V}[\dd{\varrho}{T}]$.
  В случае ограниченности $\operatorname{supp}\varrho$ это показано, например в~\cite{CarPl2000}.
  Зададим для всех  $n\in\mathbb{N}$  плотности $\varrho_n$ правилом: $\varrho_n(t)\rav\frac{n+1}{n}\varrho(t)$ если $\int_0^t\varrho(s)\,ds\leq\frac{n}{n+1}$, и
   $\varrho_n(t)\rav 0$  в противном случае. Для них при всех $T>0$ из
     $\int_T^\infty\varrho_n(t)\,dt<1$ следует
   $\bar{V}[{\cp}[\varrho_n]\equiv \bar{V}[\dd{\varrho_n}{T}]$.
 Поскольку в силу \eqref{o0},\eqref{o1}  имеет место равномерная (по $z\in\mathbb{K}$ и $\omega\in\Omega$) сходимость (при $n\to\infty$)   платежей $\cp[\varrho_n]$  к  $\cp[\varrho]$, цен $\mathcal{V}[\varrho_n]$ к $\mathcal{V}[\varrho]$,
   платежей $\dd{\varrho_n}{T}$ к $\dd{\varrho}{T}$, и
   цен  $\bar{V}[\dd{\varrho_n}{T}]$ к $\bar{V}[\dd{\varrho}{T}]$, то, переходя к пределу, для всех $\varrho\in\mathfrak{D}, T>0$
   при $\int_T^\infty\varrho_n(t)\,dt<1$    получаем
       $$
       \bar{V}[{\cp}[\varrho]]\equiv \bar{V}[\dd{\varrho}{T}].$$

       Пусть  для ${\mathcal{V}}$, заданного правилом \eqref{1588}, для некоторой
       функции  $U_*:\Omega'\to [0,1]$ выполнен хотя бы один из пунктов $(c)$, $(u)$, $(e)$, $(p)$, $(q)$. Тогда, для асимптотики $\bar{U}_*\in\mathfrak{U}'$, заданной правилом
        $\bar{U}_*(\omega)=U_*(\omega,a,b)$ $(\forall (\omega,a,b)\in\Omega')$, выполнен хотя бы один пункт
        теоремы \ref{4} с игровым отображением $\mathbb{V}=\bar{V}$, то есть выполнены все эти пункты. Делая обратную подстановку, получаем требуемое. \bo

 Эквивалентность $(u)\Leftrightarrow (e)$ для дифференциальных игр была доказана сначала только в эргодическом случае \cite{AlvBarditrue}, затем в предположении ``существования седловой точки в маленькой игре'' (условие Айзекса) \cite{MTIP}, снять это условие удалось в \cite{arxiv2016}.  Асимптотики, для распределений, отличных от равномерного и экспоненциального, ранее исследовались  в \cite{SDCP2015},
 на основе полученных там результатов недавно \cite{optima2017} для дифференциальных игр удалось показать эквивалентность условий $(u)$, $(e)$, $(c)$.

 Существование равномерных пределов цен для дифференциальных игр показано лишь  в отдельных случаях, отметим \cite{BardiGame,CQ2015,Cardal2010}.

\subsection{Стохастическая постановка}

Пусть вновь, как и в параграфе \ref{general}, заданы множества $\Omega$ и
$\mathbb{K}$, а также отображение $g:\Omega\to[0,1]$. Пусть нашлась
такая $\sigma$-алгебра ${A}$ над $\mathbb{K}$, что  отображения
$\mathbb{K}\ni z\mapsto \int_{0}^\infty \varrho(t)g(z(t))\,dt$
 ${A}$-измеримы для всех плотностей $\varrho\in\mathfrak{D}$.

Пусть для всякого $\omega\in\Omega$ заданы непустые множества  $\mathcal{L}(\omega)$ и $\mathcal{M}(\omega)$ правил для первого и второго игроков соответственно.
Теперь для всех  $\omega\in\Omega$ каждой паре $(l,m)\in\mathcal{L}(\omega)\times\mathcal{M}(\omega)$ правил игроков сопоставим некоторое
вероятностное распределение (над $(\mathbb{K},{A})$) и соответствующее ему математическое ожидание $\mathbb{M}^\omega_{lm}$.
Для всех  $\varrho\in\mathfrak{D}$ положим
\begin{eqnarray}
\mathcal{V}[\varrho](\omega)\rav\sup_{l\in \mathcal{L}(\omega)}\inf_{m\in \mathcal{M}(\omega)} \mathbb{M}^\omega_{lm} \int_{0}^\infty \varrho(t)g(z(t))\,dt\qquad \forall\omega\in\Omega.
\label{295_}
\end{eqnarray}

\begin{theorem}\label{normal11}
     Предположим, что для всех $\varrho\in\mathfrak{D},\omega\in\Omega$ были ${A}$-измеримы
     отображения
     $\mathbb{K}\ni z\mapsto \cp(z)$,
     $\mathbb{K}\ni z\mapsto \mathcal{V}[\varrho](z(T))$, и при этом, или
        \begin{eqnarray}
\mathcal{V}[\varrho](\omega)=
    {\displaystyle\sup_{l\in \mathcal{L}(\omega)}\inf_{m\in \mathcal{M}(\omega)}}
    \mathbb{M}^\omega_{lm} \Big[\int_{0}^T\!\! \varrho(t)g(z(t))\,dt\!\!\!\!\!\label{631}\\
    +
    \sup_{l'\in \mathcal{L}(z(T))}\inf_{m'\in \mathcal{M}(z(T))}\!\!\!\!\!\!
    &{\displaystyle\mathbb{M}^{z(T)}_{l'm'}}&\!\!\!\!\!\!\int_{T}^\infty\!\!\varrho(t+T)g(z_1(t))\,dt\Big]
    \nonumber
    \end{eqnarray}
имело место для всех положительных $T$, или было выполнено \eqref{579}, а  \eqref{631} имело место для всех натуральных $T$.

    Тогда, если $\mathcal{V}$ задано правилом \eqref{295_}, то для всякой функции
    $U_*:\Omega\to[0,1]$, результат теоремы~\ref{4} имеет место, в частности
     условия
    $(c)$, $(u)$, $(e)$, $(p)$, $(q)$ эквивалентны.
\end{theorem}
\underline{Доказательство теоремы~\ref{normal11}.}
Рассмотрим множество $\mathfrak{C}_0$ всех ${A}$-измери\-мых ограниченных отображений из $\mathbb{K}$ в $\mathbb{R}$.
Определим для всякого платежа $c\in\mathfrak{C}_0$
\begin{eqnarray}
V[c](\omega)\rav\sup_{l\in \mathcal{L}(\omega)}\inf_{m\in \mathcal{M}(\omega)} \mathbb{M}^\omega_{lm} c \qquad \forall\omega\in\Omega.
\label{295_1}
\end{eqnarray}
По условию, $\cp[\varrho],\dd{\varrho}{T}\in\mathfrak{C}_0$ и $V\big[\cp[\varrho]\big]=V\big[\dd{\varrho}{T}\big]$ для всех $\varrho\in\mathfrak{D}, T>0.$

Пусть $\mathfrak{C}_1$\ ---\  множество всех таких $c_1\in\mathfrak{C}_0$, что
все отображения $\mathbb{K}\ni z\mapsto V[c_1](z(T))$  ${A}$-измеримы для всех $T>0$.
По условию теоремы, $\cp[\varrho]\in\mathfrak{C}_1$ для всех плотностей $\varrho\in\mathfrak{D}.$
Тогда, из $V\big[\cp[\varrho]\big]=V\big[\dd{\varrho}{T}\big]$ следует, что
$\dd{\varrho}{T}\in\mathfrak{C}_1$
 для всех $\varrho\in\mathfrak{D}, T>0.$
Легко также проверить, что всякое постоянное отображение $c:\mathbb{K}\to \mathbb{R}$
также лежит в $\mathfrak{C}_1.$

Поскольку произвольный платеж $c\in\mathfrak{C}$  ограничен сверху некоторой константой, а значит и некоторым платежом  $c_1\in\mathfrak{C}_1$, корректно ввести
$$
 \mathbb{V}[c](\omega)\rav \inf\{{V}[c_1](\omega)\,|\,c_1\in\mathfrak{C}_1, c_1\geq c\}
 \qquad \forall \omega\in{\Omega},c\in\mathfrak{C}.
$$
 Легко проверить, см. например \cite[Lemma~1]{arxiv2016}, что это игровое отображение, при этом
 $\mathbb{V}|_{\mathfrak{C}_1}={V}|_{\mathfrak{C}_1}$. Тогда
 $\mathbb{V}\big[\cp[\varrho]\big]=\mathbb{V}\big[\dd{\varrho}{T}\big]$
 для всех $\varrho\in\mathfrak{D}, T>0.$ Теперь, из \eqref{631}, для отображения $\mathbb{V}$
 или выполнен  принцип динамического программирования на $\mathfrak{D},$
 или имеют место и  дискретный вариант этого принципа, и условие  \eqref{579}.
 Применяя в первом случае теорему~\ref{4}, а во втором\ ---\ теорему~\ref{44}, получаем требуемое. \bo

 Показанная здесь теорема применима для стохастических постановок как  в дискретном, так и непрерывном времени, и в части $(e)\Leftrightarrow (u)$ была показана в \cite{arxiv2016}. Случай непрерывного времени изучен относительно мало, отметим лишь  работы \cite{BGQ2013,mexic}. Стохастические постановки в дискретном времени, напротив, прекрасно  исследованы:
 для стохастических игр с конечными множествами состояний и действий существование предела в \eqref{1113e}
 было показано в \cite{BK}, эквивалентность $(e)\Leftrightarrow (u)$  для таких игр была доказана в \cite{MN1981}, неплохой обзор имеющихся к 2011 году результатов смотрите в \cite{Sorin2011}. Позже
 удалось показать, что уже в случае компактного множества действий  пределы в \eqref{1113u},\eqref{1113e}  могут не существовать \cite{Ziliotto2013}; также в  \cite{Ziliotto2013} был впервые предложен  метод доказательства тауберовых теорем для стохастических постановок в дискретном времени на основе построения асимптотик неподвижных точек нерасширяющихся операторов (итераций оператора Шепли). Далее, в работе \cite{Ziliotto2016},
 на основе того же подхода, для стохастических игр был показан ряд тауберовых теорем с плотностями общего вида.

\section{Обсуждение условий}
\label{prooff}

 Хотя, как показывает предыдущий параграф, теорему \ref{4} можно легко переформулировать для различных игровых постановок,
 доказывать ее мы будем вместе с еще одним утверждением, теоремой \ref{4_4}, весьма громоздкая формулировка которой, тем не менее, и упростит доказательство теоремы \ref{44},  и будет существенно удобнее для самой технически сложной части доказательства теоремы \ref{4} (см. ниже предложение~\ref{tochi}). Но сначала ввведем два обозначения.

Для любых открытого вправо интервала $[a,b)\subset\mathbb{R}$ и функции $y:[a,b)\to\mathbb{R}$, через $V_{a}^{b}[y]$ обозначим полную вариацию функции~$y$ на промежутке $[a,b)$, то есть
$$
V_{a}^{b}[y]=\sup_{k\in\mathbb{N},\ a\leq t_0<t_1<\dots<\leq t_k< b} \sum_{i=1}^k \big|y(t_i)-y(t_{i-1})\big|.
$$

Для всякой плотности $\varrho\in\mathfrak{D}$ и числа $r\in (0,1)$ введем квантиль $q[\varrho](r)$ как минимальное положительное число, для которого имеет место
$$\int_0^{q[\varrho](r)}\varrho(t)\,dt=r.$$

\begin{theorem}
    \label{4_4}
    Пусть игровое отображение~$\mathbb{V}$ таково, что
    всякому положительному $\varkappa$ найдется положительное $\gamma$ со свойством:
    для всех  $\nu\in\mathfrak{D}$, $T>1$ из
   $V_0^\infty [\nu]<\gamma, \int_0^T\nu(t)\,dt<1$ следует, что
    \begin{eqnarray}
       \label{999}
    \Big|\mathbb{V}[\dd{\nu}{T}](\omega)-\mathbb{V}\big[\cp[\nu]\big](\omega)\Big|<\varkappa  \quad\forall\omega\in\Omega.
    \end{eqnarray}
    Пусть отображение $\mathcal{V}$ определено правилом \eqref{superp}, а функция
    $U_*:\Omega\to [0,1]$ взята произвольно.

    Тогда,   каждое из условий  $(c)$, $(u)$, $(e)$, $(p)$, $(q)$ также
     эквивалентно любому из условий:
    \begin{description}
    \item[$(v)$]\
    Для всякого
    семейства плотностей $\mu_{\lambda}\in\mathfrak{D},\lambda\in(0,\lambda_*)$ со свойствами
    \begin{eqnarray}\label{lto0}
    \lim_{\lambda\to +0} \sup_{t>0}\mu_\lambda(t)&=&0,\\
    \label{lto1}
    \sup_{\lambda\in(0,\lambda_*)} V_{0}^{q[{\mu_{\lambda}}](1-\varepsilon)}[{\mu_{\lambda}}]\cdot q[\mu_{\lambda}](1-\varepsilon)&<&+\infty \quad \forall\varepsilon\in (0,1)
\end{eqnarray}
   цены $\mathcal{V}[\mu_{\lambda}]$  сходятся к~$U_*$ равномерно по $\omega\in\Omega$ при $\lambda\to+0$, то есть
\begin{equation}\label{1113g}
\lim_{\lambda\to+0} \sup_{\omega\in\Omega}
\Big|\mathcal{V}[\mu_{\lambda}](\omega)-U_*(\omega)\Big|=0.
\end{equation}
    \item[$(\exists)$]\
    Некоторое семейство плотностей $\varrho_\lambda\in\mathfrak{D}$, $\lambda\in(0,\lambda_*)$ удовлетворяет следующим трем свойствам:\\
    1) для всех  $\lambda\in(0,\lambda_*)$ имеет место
    \begin{eqnarray}
    \label{0}
    \lambda=\varrho_\lambda(0)\geq \varrho_\lambda(t)\qquad\forall t\geq 0;
    \end{eqnarray}
    2)
    для всякого   $\varepsilon>0$,
    найдутся такие положительные числа $\delta_\varepsilon<1$ и
    $\lambda_\varepsilon\leq\lambda_*$, что для всех положительных
    $\lambda<\lambda_\varepsilon,T\leq\delta_\varepsilon/\lambda$ выполнено
    \begin{eqnarray}\label{55}
    \varrho_\lambda(T)
    \geq \lambda(1-\varepsilon);
    \end{eqnarray}
    3)
     существует такое число  $r_0\in(0,1)$, что
    \begin{eqnarray}
\nonumber
    0&=&\lim_{\lambda\to+0}\sup_{\omega\in\Omega}
    \Big|\mathcal{V}[\varrho_\lambda](\omega)-U_*(\omega)\Big|\\
    &=&\lim_{\lambda\to+0}\sup_{T\in(0,q[\varrho_\lambda](r_0)),\omega\in\Omega}
    \Big|\mathcal{V}[
    (\varrho_\lambda)^{T}_{\textrm{shift}}](\omega)-U_*(\omega)\Big|.
        \label{3}
    \end{eqnarray}
    \end{description}
    \end{theorem}
    Доказательство этой теоремы (равно как теоремы \ref{4}) начнется в параграфе~\ref{proofg}, при этом нам потребуется также доказать (см. параграф~\ref{doctochi})
\begin{proclaim}
    \label{tochi}
    Пусть выполнены предположения теоремы \ref{4_4}, и для некоторой функции $U_*\in\mathfrak{U}$ имеет место условие $(\exists)$.

    Тогда для любых положительных $\varepsilon,M$ найдется такое положительное число  $\gamma>0$, что для всех плотностей $\mu\in\mathfrak{D}$ из
    \begin{eqnarray}\label{w1}
    \sup_{t>0}\mu(t)<\gamma,\quad V_{0}^{q[{\mu}](1-\varepsilon)}[{\mu}]\cdot q[\mu](1-\varepsilon)\leq M,
    \end{eqnarray}
    следует
    $\mathcal{V}[\mu](\omega)\geq U_*(\omega)-12\varepsilon$  для всех $\omega\in\Omega$.
\end{proclaim}

Вернемся к обсуждению теорем.

 Теоремы~\ref{4},\ref{44},\ref{4_4} позволяют, проверив наличие равномерной асимптотики $U_*$ для цен,
 соответствующих какому-то заданному, тестирующему семейству плотностей $\varrho_{\lambda}$ (например
 \eqref{1113u}--\eqref{1113m}), гарантировать ту же асимптотику цен, то есть
 \eqref{1113g}, для семейства плотностей $\mu_{\lambda}$, удовлетворяющего
 \eqref{lto0},\eqref{lto1}.  Ранее в литературе  в качестве  тестирующих рассматривались лишь  семейства
 экспоненциальных и(или) равномерных распределений.  Условие $(\exists)$ описывает более мягкие
 требования, предъявляемые к тестирующему семейству, в частности, как следует из условия $(p)$,  достаточно проверить сдвиги по времени произвольной степенной плотности (см. \eqref{1113m}).

 Не менее интересны достаточные для \eqref{1113g} требования на   семейство плотностей $\mu_{\lambda},\lambda\in(0,\lambda_*)$, требования, гарантирующие для него  ту же  асимптотику цен, что и для тестирующего семейства.
  В целом ряде постановок для обеспечения \eqref{1113g}  используются существенно более слабые, нежели \eqref{lto0}--\eqref{lto1} предположения.
   В частности, в \cite{Ziliotto2016} показано, что для любой стохастической игры с конечным числом  состояний и действий,  от семейства плотностей $\mu_{\lambda}$
    достаточно потребовать,   хоть при каком-то  положительном~$p$, выполнения
  \begin{eqnarray}\label{lto3}
   \lim_{\lambda\to+0}V_0^\infty [(\mu_{\lambda})^p]=0.
    \end{eqnarray}
    Условие \eqref{lto3} с $p=1$ достаточно  в детерминированных процессах с дискретным временем \cite{Renault2013}, для марковских процессов с компактным множеством состояний~\cite{RenaultVenel2012}. Для задач управления в эргодическом случае в \cite{LiQixRenault}
    предложено сходное (см. \cite[Proposition~3.6]{LiQixRenault}) условие \eqref{501}.

    Однако в \cite{Ziliotto2017} построен пример управляемой марковской цепи с таким семейством монотонно невозрастающих плотностей $\mu_{\lambda}$, длч которого выполнено требование
    \eqref{lto0} (а значит и \eqref{lto3} для всякого $p>0$), имеет место сходимость \eqref{1113u} (условие $(u)$), но сходимость
    \eqref{1113g} (условие $(c)$), тем не менее, места не имеет.
    Отсюда следует, что  в условиях теорем~\ref{44} и \ref{4_4}, требование \eqref{lto1} не может быть заменено на \eqref{lto3} даже в случае одного игрока.

    Основное требование, предъявляемое в   теоремах~\ref{4},\ref{44},\ref{4_4} к самой игре --- та или иная форма принципа динамического программирования. Без какой-либо вариации этого принципа       уже результат следствия~\ref{u1} может не иметь место, см. соответствующий контрпример  в \cite{MatSb}. Впрочем, предложенная здесь (см. \eqref{999}) асимптотическая версия  принципа динамического программирования по-видимому может быть ослаблена, например, в \cite{arxiv2016} требовался лишь её  дискретный вариант, и  только для экспоненциального и равномерного распределений.
    Более того,  применяя принципы суб- и супероптимальности \cite[Subsection~3.3]{arxiv2016},    (знаки ``$\geq$'' или``$\leq$'' вместо ``$=$'' в  \eqref{435})      можно показать односторонние тауберовы теоремы \cite[Proposition~3-6]{arxiv2016},

    В  теоремах~\ref{4},\ref{44},\ref{4_4} также требуется, чтобы игровое отображение $\mathbb{V}$ было задано на множестве  всех ограниченных функций из $\mathbb{K}$ в $\mathbb{R}.$
     Это не всегда удобно например в стохастическом случае, но это условие всегда можно обеспечить, доопределив игровое отображение (см. доказательство теоремы \ref{normal11},\
     а также  \cite[Lemma~1]{arxiv2016}).

    Также хорошо известно,  что уже в самых простых управляемых системах, тауберова теорема не верна в случае  применения поточечного (а не равномерного) предела цен; см.  пример в~\cite{barton}.

    Перейдем к доказательству указанных выше теорем.

\section{Доказательство теорем~\ref{4}, \ref{44} и \ref{4_4}}
\label{proofg}

\subsection{Доказательство теорем~\ref{4} и \ref{4_4}}
Мы докажем эти теоремы одновременно, по следующей схеме:
 $(u)\Rightarrow (\exists),$ $(e)\Rightarrow (\exists),$ $(c)\Rightarrow (u),(c)\Rightarrow (e),$
 $(\exists)\Rightarrow (v)\Rightarrow (c)\Rightarrow (q)\Rightarrow (p)\Rightarrow (\exists).$
Начнем с
\subsection*{Доказательство $(u)\Rightarrow (\exists),$ $(e)\Rightarrow (\exists)$, $(c)\Rightarrow (u),$
    $(c)\Rightarrow (e)$}
Напомним, что   для всех $t\geq 0,\lambda>0$
$\varpi_{1/\lambda}(t)=\lambda 1_{[0,1/\lambda]}(t)$ и $\pi_\lambda(t)=\lambda e^{-\lambda t}$.
Легко видеть, что условие \eqref{0}
 выполнено для каждого из этих семейств.

 Примем  $\lambda_*=1/2.$ Каждому положительному $\varepsilon<\lambda_*$ сопоставим $\lambda_\varepsilon\rav\varepsilon,\delta_\varepsilon\rav\varepsilon$;
 теперь для всех положительных $\lambda<\lambda_\varepsilon$ и
$T\leq \delta_\varepsilon/\lambda$ выполнено
  $1=\varpi_{1/\lambda}(T)/\lambda\geq\pi_{\lambda}(T)/\lambda\geq \pi_\lambda(\delta_\varepsilon/\lambda)/\lambda=e^{-\delta_\varepsilon}>1-\varepsilon,$
  то есть \eqref{55}.

 Примем  $r_0=e^{-1/2}<1/2.$ Поскольку
 $(2\lambda)^{-1}=q[\pi_{\lambda}](e^{-1/2})=q[\varpi_{1/\lambda}](1/2)$ для всех $\lambda>0$,  то   для всех положительных $T<(2\lambda)^{-1}$ имеет место
$$
(\varpi_{1/\lambda})^{T}_{\textrm{shift}}=\varpi_{\frac{1}{\lambda}-T},
 \quad
  (\pi_{\lambda})^{T}_{\textrm{shift}} =\pi_{\lambda}.
  $$
 Теперь  пределы в~\eqref{3} существуют, равны, и равномерны (по $\omega\in\Omega$ и по  $\omega\in\Omega$, $T\in[0,q[\varrho_\lambda](r_0)]$) или при $\varrho_\lambda=\varpi_{1/\lambda}$ или при $\varrho_\lambda=\pi_\lambda$
в зависимости от того какое из свойств, $(u)$ или~$(e)$, выполнено.
Таким образом,   $(u)\Rightarrow (\exists),$ $(e)\Rightarrow (\exists)$ доказаны.

 Для доказательства импликаций $(c)\Rightarrow (u),$
 $(c)\Rightarrow (e)$
 заметим, что плотности $\pi_1,\varpi_1$\ ---\ кусочно-непрерывны, далее,
 для всякого $\lambda>0$ выполнено
 $$
 \big(\pi_1\big)^\lambda_{\textrm{scale}}=\pi_\lambda,\quad
 \big(\varpi_1\big)^\lambda_{\textrm{scale}}=\varpi_{1/\lambda}
 \qquad \forall \lambda>0.
 $$
 Теперь, подставляя в \eqref{1113b} $\mu=\pi_1$ и $\mu=\varpi_1$, имеем из $(c)$ как $(e)$, так и $(u)$.

\subsection*{Доказательство  $(\exists)\Rightarrow (v)$}
Пусть плотности $\varrho_\lambda$ удовлетворяют условию $(\exists)$, а для
плотностей~$\mu_\lambda,\lambda\in(0,\lambda_*)$ выполнены \eqref{lto0},\eqref{lto1}.

Каждому положительному~$\varepsilon$ сопоставим число
$$M_\varepsilon=\sup_{\lambda\in(0,\lambda_*)} V_{0}^{q[{\mu_{\lambda}}](1-\varepsilon)}[{\mu_{\lambda}}]\cdot q[\mu_{\lambda}](1-\varepsilon)\in\mathbb{R}.$$
По предложению~\ref{tochi}, можно найти такое $\gamma^+_\varepsilon>0$, что для всех
$\lambda\in(0,\lambda_*)$ из
$\sup_{t>0}\mu_\lambda (t)<\gamma^+_\varepsilon$
следует
 $U_*(\omega)-12\varepsilon\leq\mathcal{V}[\mu_\lambda](\omega)$ для всех $\omega\in\Omega$.

Рассмотрим новую функцию мгновенной полезности $g^-\rav 1-g$  и еще одно отображение $V^-[c]\rav-\mathbb{V}[-c].$
Легко проверить, что $V^-$ --- игровое отображение, и  для всех плотностей~$\varrho\in\mathfrak{D}$,
 \begin{eqnarray*}
 1-\mathcal{V}[\varrho](\omega)=V^-[1-\cp[\varrho]](\omega)\qquad \forall\omega\in\Omega.
\end{eqnarray*}
Условия предложения~\ref{tochi} выполнены для игрового отображения~$V^-$ с семейством плотностей $\varrho_\lambda$, функцией мгновенной полезности~$g^-$ и асимптотикой $1-U_*$.
Тогда найдется такое
 $\gamma^-_\varepsilon>0$, что для всех
$\lambda\in(0,\lambda_*)$ из
$\sup_{t>0}\mu_\lambda (t)<\gamma^-_\varepsilon$
следует
 $1-\mathcal{V}[\mu_\lambda](\omega)\geq 1-U_*(\omega)-12\varepsilon$ для всех $\omega\in\Omega$.

Итак, для для всех
$\lambda\in(0,\lambda_*)$ из
$\sup_{t>0}\mu_\lambda (t)<\min(\gamma^+_\varepsilon,\gamma^-_\varepsilon)$ следует
 $U_*(\omega)-12\varepsilon\leq \mathcal{V}[\mu_\lambda](\omega)\leq U_*(\omega)+12\varepsilon$.
 Теперь
 условие \eqref{1113g} следует из \eqref{lto0}.

\subsection*{Доказательство $(v) \Rightarrow (c)$}
Рассмотрим некоторую кусочно-непрерывную на всяком компакте плотность $\bar\varrho\in\mathfrak{D}$.

Вместе с некоторым  положительным $\varepsilon$ зафиксируем некоторое, достаточно большое натуральное число $n>3$, для которого ${10}/{n}<\varepsilon$.

Зададим числа
 \begin{eqnarray*}
 r_n\rav q[\varrho](1/n),\quad s_n\rav q[\varrho](1-1/n).
 \end{eqnarray*}
Поскольку~$\varrho$ кусочно-непрерывна на $[r_n,s_n]$, найдется кусочно-постоянная функция $\mu_n:\mathbb{R}_+\to \mathbb{R}_+$, $\operatorname{supp}\mu_n\subset{[r_n,s_n]}$,
для которой $\int_{r_n}^{s_n}\mu_n(t)\,dt=\int_{r_n}^{s_n}\varrho(t)\,dt=\frac{n-2}{n}$,
 $\int_{r_n}^{s_n}|\mu_n(t)-\varrho(t)|\,dt<1/n$.
В частности, конечна ее полная вариация на отрезке $[r_n,s_n]$.
В силу $\int_0^\infty \mu_n(t)\,dt=\frac{n-2}{n}$,
установим
 \begin{eqnarray*}
 \bar{\mu}\rav\frac{n}{n-2}\mu_n\in\mathfrak{D},\quad
 M= s_n V_{0}^{\infty}[\bar{\mu}]=\frac{n}{n-2}s_n V_{0}^{\infty}[{\mu}_n]\in\mathbb{R}.
 \end{eqnarray*}
Сейчас для всех $\lambda>0$ мы имеем
 \begin{eqnarray}
    V_{0}^{\infty}[\bar{\mu}^{\lambda}_{\textrm{scale}}]\cdot q[\bar{\mu}^{\lambda}_{\textrm{scale}}](1-\varepsilon)&=&
\lambda V_{0}^{\infty}[\bar{\mu}]\cdot\frac{q[\bar{\mu}](1-\varepsilon)}{\lambda}\leq M,
\label{8111}\\
  \int_{0}^\infty \Big|\varrho^\lambda_{\textrm{scale}}(t)-\bar{\mu}^{\lambda}_{\textrm{scale}}(t)\Big|\,dt&=&
  \int_{0}^\infty \Big|\varrho(t)-\bar{\mu}(t)\Big|\,dt\nonumber\\
  &<&\frac{3}{n}+
  \int_{r_n}^{s_n} \Big(\frac{n}{n-2}-1\Big)\varrho(t)\,dt=\frac{5}{n}<\frac{\varepsilon}{2}.\nonumber
     \end{eqnarray}
Таким образом,  в силу \eqref{o1}, показано
 \begin{eqnarray}
  \sup_{\omega\in\Omega}\big|\mathcal{V}[\varrho^\lambda_{\textrm{scale}}](\omega)-\mathcal{V}[\bar{\mu}^{\lambda}_{\textrm{scale}}](\omega)\big|\leq\frac{\varepsilon}{2},\label{811}
       \end{eqnarray}
более того, в силу \eqref{8111} для семейства
$\bar{\mu}^{\lambda}_{\textrm{scale}},\lambda_n>0$
	выполнено \eqref{lto1}.
	
Для всех положительных $\lambda$ имеем
$$\sup_{t>0}\bar{\mu}^{\lambda}_{\textrm{scale}}(t)\leq
V_{0}^{\infty}[\bar{\mu}^{\lambda}_{\textrm{scale}}]=\lambda V_{0}^{\infty}[\bar\mu_n],$$
таким образом условие \eqref{lto0} для семейства
$\bar{\mu}^{\lambda}_{\textrm{scale}},\lambda>0$
	также выполнено.
Теперь, для этого семейства, в силу~$(v)$, показано  \eqref{1113m}, в частности для некоторого положительного $\lambda_n$ имеет место
   $$\sup_{\omega\in\Omega}|\mathcal{V}[\bar{\mu}^{\lambda}_{\textrm{scale}}](\omega)-U_*(\omega)|<\frac{\varepsilon}{2}$$
для всех положительных $\lambda<\lambda_n.$
Для таких~$\lambda$, с учетом~\eqref{811} получаем
   $$\sup_{\omega\in\Omega}|\mathcal{V}[\varrho^\lambda_{\textrm{scale}}](\omega)-U_*(\omega)|<\varepsilon.$$
   В силу произвольности выбора положительного~$\varepsilon$, импликация  $(v) \Rightarrow (c)$ доказана.

\subsection*{Доказательство $(c)\Rightarrow (q)\Rightarrow (p)\Rightarrow (\exists)$}
  Покажем сначала, что для всякой степенной плотности $\mu\in\mathfrak{D}$ выполнено
   \begin{eqnarray}\label{glav}
    \limsup_{\lambda\to+0} \sup_{\omega\in\Omega}
    \Big|\mathcal{V}[\mu^\lambda_{\textrm{scale}}](\omega)-U_*(\omega)\Big|=
    \limsup_{T\to +\infty} \sup_{\omega\in\Omega}
    \Big|\mathcal{V}[\mu^T_{\textrm{shift}}](\omega)-U_*(\omega)\Big|.
   \end{eqnarray}
    Действительно, пусть задана некоторая степенная плотность $\mu\in\mathfrak{D},$ то есть для некоторых $\alpha,\beta>0,\gamma>1$ имеет место $\mu(t)=(\alpha+\beta t)^{-\gamma}$ для всех $t\geq 0.$
 Теперь из $\mu\in\mathfrak{D}$ следует, что
   \begin{eqnarray*}
  \int_{T}^\infty \mu(t)\,dt=\frac{\int_{T}^\infty \mu(t)\,dt}{\int_{0}^\infty \mu(t)\,dt}=\frac{\alpha^{\gamma-1}}{(\alpha+\beta T)^{\gamma-1}}\qquad\forall T>0.
  \end{eqnarray*}

Всякому положительному  $\lambda<1$ сопоставим $T(\lambda)=\frac{\alpha}{\beta\lambda}(1-\lambda)$;
при этом $T(+0)=+\infty$ и найдется обратное отображение $T\mapsto \lambda(T)$, действующее по правилу
 $(0,\infty)\ni T\mapsto\lambda(T)=\frac{\alpha}{\alpha+\beta T}\in(0,1)$.
Теперь,
   \begin{eqnarray}
\mu^{\lambda(T)}_{\textrm{scale}}(t)&=&\frac{\lambda(T)}{(\alpha+\beta \lambda(T) t)^{\gamma}}=
\frac{\lambda^{1-\gamma}(T)}{(\frac{\alpha}{\lambda(T)}+\beta t)^{\gamma}}\nonumber\\
&=&
\frac{\alpha^{1-\gamma}}{(\alpha+\beta T)^{1-\gamma}(\alpha+\beta T+\beta t)^{\gamma}}\nonumber\\&=&
\frac{\mu(t+T)}{\int_{T}^\infty \mu(t)\,dt}=
  \mu^{T}_{\textrm{shift}}(t)\qquad \forall t\geq 0.\label{glav1}
  \end{eqnarray}
  Поскольку $T(+0)=+\infty$, то \eqref{glav} показано для всех степенных плотностей.

Покажем  $(c)\Rightarrow (q)$: действительно,
  условие~$(c)$ гарантирует, что
 $\mathcal{V}[\mu^\lambda_{\textrm{scale}}]$ сходятся равномерно (по $\Omega$) к  $U_*$ при $\lambda\to+0,$
 теперь из \eqref{glav} следует, что
 $\mathcal{V}[\mu^T_{\textrm{shift}}]=\mathcal{V}[\mu^{\lambda(T)}_{\textrm{scale}}]$ равномерно сходятся к~$U_*$ при $T\to+\infty.$ Таким образом, импликация $(c)\Rightarrow (q)$ показана.

Импликация  $(q)\Rightarrow (p)$ непосредственно следует из существования хотя бы одной степенной плотности, а ее можно  задать например  равенством $\mu(t)=(1+t)^{-2}$ для всех $t\geq 0$.

Докажем импликацию $(p)\Rightarrow (\exists)$.

    Пусть задана некоторая степенная плотность $\mu\in\mathfrak{D}.$ то есть для некоторых $\alpha,\beta,\gamma>0$ имеет место $\mu(t)=(\alpha+\beta t)^{-\gamma}$ для всех $t\geq 0.$

    Для всех $\lambda>0$ примем
    $\varrho_\lambda=\mu^{\lambda/\mu(0)}_{\textrm{scale}}$.
    В силу  \eqref{glav1} для всех  положительных $\lambda<\mu(0)$
    $$\varrho_\lambda=\mu^{\lambda/\mu(0)}_{\textrm{scale}}
    =\mu^{T(\lambda/\mu(0))}_{\textrm{shift}}
        =\mu^{\frac{\alpha}{\beta\lambda}(\mu(0)-\lambda)}_{\textrm{shift}}.$$
    Поскольку при этом  $\mu(0)/\lambda$ стремится к $+\infty$ при $\lambda\to+0$, то из \eqref{1113m} следует
     \begin{equation}\label{1113mm}
    \lim_{\lambda\to+0} \sup_{\omega\in\Omega}
    \Big|\mathcal{V}[\varrho_\lambda](\omega)-U_*(\omega)\Big|=0.
    \end{equation}

    Легко видеть, что
    $\varrho_\lambda(0)=\mu^{\lambda/\mu(0)}_{\textrm{scale}}(0)=\lambda/\mu(0)\cdot\mu(0)=\lambda,$
    и для семейства плотностей $\varrho_\lambda$ условие \eqref{0}
    выполнено при $\lambda_*=\mu(0)$.

    Далее, в силу непрерывности плотности $\mu$, для каждого положительного $\varepsilon<1$ найдется такое положительное $\delta_\varepsilon$, для которого $\mu\big(t/\mu(0)\big)/\mu(0)=1-\varepsilon.$
    Теперь для всех положительных $\lambda<1$, $T<\delta_\varepsilon/\lambda$ получаем
    $$ \frac{\varrho_\lambda(T)}{\varrho_\lambda(0)}=
    \frac{\mu^{\lambda/\mu(0)}_{\textrm{scale}}(T)}{\mu^{\lambda/\mu(0)}_{\textrm{scale}}(0)}=
    \frac{\mu(T\lambda/\mu(0))}{\mu(0)}>\frac{\mu(\delta_\varepsilon/\mu(0))}{\mu(0)}=
    1-\varepsilon,
    $$
    Итак, \eqref{55} проверено.

    Рассмотрим теперь для некоторых положительных $T$ и $\lambda'$ $(\lambda'<\mu(0))$ плотность
    $(\varrho_{\lambda'})^{T}_{\textrm{shift}}$.  Поскольку $\varrho_{\lambda'}$ также является степенной плотностью, к ней применимо \eqref{glav1}, откуда получаем
    $$
     \left(\varrho_{\lambda'}\right)^{T}_{\textrm{shift}}=
      \left(\mu^{\lambda'/\mu(0)}_{\textrm{scale}}\right)^{T}_{\textrm{shift}}=
      \left(\mu^{\lambda'/\mu(0)}_{\textrm{scale}}\right)^{\lambda(T)}_{\textrm{scale}}=
      \mu^{\lambda(T)\lambda'/\mu(0)}_{\textrm{scale}}=\varrho_{\lambda(T)\lambda'}=
      \varrho_{\frac{\alpha\lambda'}{\alpha+\beta T}}.
    $$
     Поскольку для всех $T>0$ $\frac{\alpha\lambda'}{\alpha+\beta T}\leq \lambda'$, из
     \eqref{1113mm} следует и \eqref{3}.
Итак, условие $(\exists)$ показано.

Доказательство теорем~\ref{4} и~\ref{4_4} завершено.

     \subsection{Доказательство теоремы~\ref{44}}

          Покажем сначала  $(u)\Leftrightarrow (u')$.
 Как доказано в \cite[(8a)]{arxiv2016}, для всех  $T>0$, $r>1$, $\omega\in\Omega$, $z\in\mathbb{K}$ выполнено
\begin{eqnarray*}
	\Big|\mathcal{V}[\varpi_{T}](\omega)-\mathcal{V}[\varpi_{rT}](\omega)\Big|\leq 2(r-1).
\end{eqnarray*}
Теперь, поскольку всякое  положительное $T$ можно представить в виде $T=n-s$ для некоторых $n\in\mathbb{N},s\in[0,1),$
отсюда следует
$$\Big|\mathcal{V}[\varpi_{T}](\omega)-\mathcal{V}[\varpi_{n}](\omega)\Big|\leq 2\Big(\frac{n}{T}-1\Big)\leq\frac{2}{T}\qquad \forall\omega\in\Omega.$$
Итак, $(u')\Leftrightarrow (u)$ доказано.

     Теперь, если мы покажем, что
      из \eqref{579} и дискретного варианта принципа динамического программирования  следует
      \eqref{999}, то теорема~\ref{44} будет напрямую следовать из теоремы~\ref{4_4}.

 Рассмотрим  произвольные плотность $\varrho\in\mathfrak{D}$ и $r\in(0,1)$, для которых
 $\varrho^r_{\textrm{shift}}$ определена.
   Докажем, что
        \begin{eqnarray}\label{999r}
        \Big|\cp[\varrho](z)-
        \int_{r}^{\infty}\varrho(t)\,dt\cdot
        \cp[\varrho^r_{\textrm{shift}}](z)\Big|
        \leq V_0^\infty [\varrho].
        \end{eqnarray}
   Введем новую плотность $\bar\varrho\in\mathfrak{D}$ правилом: для всех $n\in\mathbb{N}$
\begin{eqnarray*}
           \bar\varrho(t)=\left\{\begin{array}{ll}
           {\displaystyle\frac{1}{r}\int_{n-1}^{n-1+r}\varrho(s)\,ds,} & \textrm{при\ }  t-n+1\in[n-1,n-1+r);\\
            {\displaystyle\frac{1}{1-r}\int_{n-1+r}^{n}\varrho(s)\,ds,}& \textrm{при \ } t\in[n-1+r,n).
            \end{array}\right.
\end{eqnarray*}
Благодаря \eqref{579},   по построению имеем
   $\cp[\bar\varrho]\equiv\cp[\varrho],\cp[\bar\varrho^r_{\textrm{shift}}]\equiv\cp[\varrho^r_{\textrm{shift}}].$
Более того,  $V_0^\infty [\bar\varrho]\leq V_0^\infty [\varrho]$.
  В силу этих неравенств достаточно доказать \eqref{999r} для так построенной $\bar\varrho$, а это можно сделать прямым подсчетом:
   для всех
   $z\in\mathbb{K}$
     \begin{eqnarray*}
     &\ &\Big|\cp[\bar\varrho](z)-\int_{r}^\infty\bar\varrho(s)\,ds\cdot
     \cp[\bar\varrho^r_{\textrm{shift}}](z)\Big|\\&=&
       \Big|\int_{0}^\infty \bar\varrho(s)z(s)\,ds-\int_{0}^\infty \bar\varrho(s+r)z(s)\,ds\Big|\\
       &\leqref{579}&
    \sum_{k=1}^\infty \Big|\int_0^r\bar\varrho(s+k)\,ds-\int_0^r\bar\varrho(s+k+1)\,ds\Big|z(k)\\
       &\leq& r
           \sum_{k=1}^\infty \Big|\bar\varrho(k)-\bar\varrho(k+1)\Big|< V_0^\infty [\bar\varrho].
        \end{eqnarray*}
       Теперь \eqref{999r} показано для всех $\varrho\in\mathfrak{D}$.

       Воспользовавшись  \eqref{o0}, подставляя в \eqref{999r} вместо~$\varrho$ некоторый ее сдвиг
       $\varrho^n_{\textrm{shift}}$, получаем для всех $n\in\mathbb{N},r\in(0,1)$, для которых $\varrho^{n+r}_{\textrm{shift}}$ существует,
               \begin{eqnarray}\label{999s}
               &\ &\Big|\int_{n}^\infty \varrho(t)\,dt\cdot \mathcal{V}[\varrho^n_{\textrm{shift}}](\omega)-
               	\int_{n+r}^\infty \varrho(t)\,dt\cdot\mathcal{V}[\varrho^{n+r}_{\textrm{shift}}](\omega)\Big|\\
               &\leq& \int_{n}^\infty \varrho(t)\,dt \cdot
                                V_0^\infty \big[\varrho^n_{\textrm{shift}}\big]=
               	V_{n}^\infty [\varrho]\leq V_{0}^\infty [\varrho]\qquad \forall \omega\in\Omega\nonumber.
               \end{eqnarray}

   Покажем, наконец, \eqref{999}.
   Cнова зафиксируем некоторую плотность $\varrho\in\mathfrak{D}$, а с ним некоторое
      число $T>1$.
     Для некоторых $n\in\mathbb{N}\cup\{0\}$ и
     $s\in[0,1)$ выполнено $T=n+s.$
     Достаточно рассмотреть случай, когда $\int_0^T \varrho(t)\, dt<1$, $\int_n^T \varrho(t)\, dt>1$.

     Примем $R=V_0^\infty [\varrho]$. Отметим, что поскольку $\varrho(t)$ стремится к нулю при $t\to\infty$, то $\sup_{t\geq 0}\varrho(t)\leq R,$ откуда
     \begin{eqnarray}
     \label{9577}
\int_0^\infty  [\varrho(t+n)-\varrho(t+n+s)]\,dt=     \int_{n}^{T}\varrho(t)\,dt\leq sR<R.
     \end{eqnarray}
     Используя дискретный вариант принципа динамического программирования,
     отметим:
 с функцией~$\mathcal{V}[\varrho]$
совпадает функция, значение отображения~$\mathbb{V}$ для платежа
     \begin{eqnarray*}
         \mathbb{K}\ni z\mapsto\int_{0}^{n}\varrho(t)g\big(z(t)\big)\,dt+
         \int_{n}^{\infty}\varrho(t)\,dt\cdot
         \mathcal{V}\big[\varrho^{n}_{\textrm{shift}}\big]\big(z(n)\big)\\
         \ravref{579}    \int_{0}^{n}\varrho(t)g\big(z(t)\big)\,dt+
         \int_{n}^{\infty}\varrho(t)\,dt\cdot
         \mathcal{V}\big[\varrho^{n}_{\textrm{shift}}\big]\big(z(T)\big).
     \end{eqnarray*}
     Теперь, из \eqref{9577}, между функциями $\mathcal{V}[\varrho]-R$ и $\mathcal{V}[\varrho]+R$ находится функция, значение
     отображения~$\mathbb{V}$ для платежа
     $$
         \mathbb{K}\ni z\mapsto  \int_{0}^{n+s}\varrho(t)g\big(z(t)\big)\,dt+
         \int_{n}^{\infty}\varrho(t)\,dt\cdot
     \mathcal{V}\big[\varrho^{n}_{\textrm{shift}}\big]\big(z(T)\big).
     $$
     Благодаря~\eqref{999s},  между $\mathcal{V}[\varrho]-2R$ и $\mathcal{V}[\varrho]+2R$ лежит  значение отображения~$\mathbb{V}$ для платежа
     $$
             \mathbb{K}\ni z\mapsto\dd{\varrho}{T}=\int_{0}^{n+s}\varrho(t)g\big(z(t)\big)\,dt+\int_{n+s}^{\infty}\varrho(t)\,dt\cdot
     \mathcal{V}\big[\varrho^{n+s}_{\textrm{shift}}\big]\big(z(T)\big),
     $$
     то есть,
     $$\Big|\mathbb{V}[\dd{\varrho}{T}](\omega)-\mathcal{V}[\varrho](\omega)\Big|\leq 2R=2V_0^\infty [\varrho]\qquad\forall \omega\in\Omega,T>1,\varrho\in\mathfrak{D}.$$

     Итак, условие~\eqref{999} показано. Теперь все следует из теоремы~\ref{4_4}.

\section{Доказательство предложения 1}
\label{doctochi}

Для всякого платежа $\mathbb{K}\ni z\mapsto c(z)\in \mathbb{R}$ введем следующее обозначение:
$$\ssp{c(z)}{\omega}\rav \mathbb{V}[c](\omega)\qquad\forall\omega\in\Omega.$$
Например, выражения
$$ \sspk{\int_{0}^h a(t)g(z(t))\,dt+ U_*(z(h))}{\omega},\
 \sspk{\int_{0}^h b(t)g(z_1(t))\,dt+ \ssp{c(z)}{z_1(h)}}{z(h')}$$
 являются значениями игрового отображения $\mathbb{V}$  для, соответственно, платежа
 $\mathbb{K}\ni z\mapsto \int_{0}^h a(t)g(z(t))\,dt+ U_*(z(h))\in\mathbb{R}$ в точке $\omega$ и
 платежа
 $$\mathbb{K}\ni z_1\mapsto \int_{0}^h b(t)g(z_1(t))\,dt+\mathbb{V}[c](z_1(h))\in\mathbb{R}$$ в точке $z(h')$.
 Отметим, что, подобно вспомогательным переменным внутри интеграла,
 произвольно выбрать  можно  и
 находящиеся внутри квадратных скобок символы переменных, если их значения пробегают всё $\mathbb{K}$; таким образом
  $\ssp{\ssp{U_*(z_5(1))}{z_{2}(1)}}{\omega}$
  автоматически совпадает  с $\ssp{\ssp{U_*(z_1(1))}{z(1)}}{\omega}$.

\subsection*{Шаг 1. Выбор констант}

        Зафиксируем некоторые положительные $M,\varepsilon$. Можно считать, что $\varepsilon<1/10, M>1.$

Выберем из условия~$(\exists)$ числа $\lambda_\varepsilon,\delta_\varepsilon$, а также $r_0\in(0,1)$.
Уменьшая при необходимости $\delta_\varepsilon,r_0,\lambda_*$ мы можем полагать, что  $\delta_\varepsilon/2=r_0<1/2,\lambda_*=\lambda_\varepsilon.$
Зафиксируем такие $\delta_\varepsilon,r_0,\lambda_*,\lambda_\varepsilon.$

По условию \eqref{55} для всех положительных $\lambda<\lambda_\varepsilon$
и $t\leq \delta_\varepsilon/\lambda$ имеет место
  $\varrho_\lambda(t)\geq \lambda(1-\varepsilon),$
откуда
  $$\int_{0}^{\delta_\varepsilon/\lambda} \varrho_\lambda(t)\,dt\geq\int_{0}^{\delta_\varepsilon/\lambda}\lambda(1-\varepsilon)\,dt=(1-\varepsilon)\delta_\varepsilon>\delta_\varepsilon/2=r_0.$$
Следовательно, $q[\varrho_\lambda](r_0)<\delta_\varepsilon/\lambda,$ и тогда,
для всех положительных $T<q[\varrho_\lambda](r_0)$ получаем
  $\varrho_\lambda(T)\geq \lambda(1-\varepsilon)$,
 \begin{eqnarray*}
  1-\int_{0}^T\varrho_\lambda(t)\,dt\leq 1-\int_{0}^T\lambda(1-\varepsilon)\,dt=1-T\lambda +T\lambda\varepsilon;
 \end{eqnarray*}
воспользовавшись неравенствами $U_*\geq 0,\ g\geq 0$,  а также
  $\varrho_\lambda(t)\leq\lambda$, мы для всех  $z\in\mathbb{K}$ гарантируем
 \begin{eqnarray*}
  \int_0^{T}\varrho_\lambda(t) g(z(t))\,dt+\int_{T}^{\infty} \varrho_\lambda(t)\,dt\cdot U_*(z(T))\\
  \leq \int_0^{T}\lambda g(z(t))\,dt+(1-T\lambda)U_*(z(T))+2T\lambda\varepsilon.
 \end{eqnarray*}
 Итак, для всех положительных $\lambda<\lambda_*$,$T\leq r_0/\lambda<q[\varrho_\lambda](r_0)$ и любого $\omega\in\Omega$ доказано
 \begin{eqnarray}
 \label{1157}
    \ssp{\int_0^{T}\varrho_\lambda(t) g(z(t))\,dt+\int_{T}^{\infty} \varrho_\lambda(t)\,dt\cdot U_*(z(T))}{\omega}\\
    \leq
  \ssp{\int_0^{T}\lambda g(z(t))\,dt+(1-T\lambda)U_*(z(T))}{\omega}+2T\lambda\varepsilon.\nonumber
 \end{eqnarray}

Выберем натуральное число
 $k>\varepsilon/r_0$, достаточно большое, чтобы обеспечить
               $$k\varepsilon>\ln\frac{1}{\varepsilon},\quad
                            k\varepsilon\ln(1+\varepsilon)>M.$$
Теперь, воспользовавшись неравенством  $1-s<e^{-s}$ при $s=\frac{1}{k^2}\ln\frac{1}{\varepsilon}$, мы также гарантировали
$k(1-\varepsilon^{1/k^2})<\frac{1}{k}\ln\frac{1}{\varepsilon}<\varepsilon$.

Зафиксируем такое~$k$.  Установим
   \begin{equation}
   \label{M}
p\rav\varepsilon^{1/k^2},\ \delta\rav 1-p<\frac{\varepsilon}{k}<r_0<1/2,\
\varkappa\rav\varepsilon(1-p).
   \end{equation}
Отметим, что эти числа  выбраны так, что   $\delta+p\delta+\dots+p^{k^2-1}\delta=1-\varepsilon,$
  $\delta+p\delta+p^2\delta+\dots=1,$ $\varkappa(1+p+p^2+\dots)=\varepsilon.$

Воспользовавшись~\eqref{3}, для некоторого положительного $\hat{\lambda}<\lambda_*$ при любом выборе положительных
     $\lambda<\hat\lambda$, $T\leq\delta/\lambda$, в силу $\delta/\lambda\leq q[\varrho_\lambda](\delta)<q[\varrho_\lambda](r_0)$, имеем
   \begin{equation}
   \label{1053}
   U_*-\mathcal{V}[{\varrho_{\lambda}}]<\varkappa,\qquad
   \mathcal{V}[(\varrho_\lambda)^{T}_{\textrm{shift}}]-U_*<\varkappa.
   \end{equation}

Поскольку предположения теоремы~\ref{4_4} выполнены, то для некоторого положительного  $\check{\lambda}<\hat{\lambda}$, для всех  $\nu\in\mathfrak{D},T>1$ из
$V_0^\infty[\nu]\leq \check{\lambda}$, $\int_0^T\nu(t)\,dt<1$ имеем для всех  $\omega\in\Omega:$
   \begin{equation}
   \label{105353}
      \bigg|
   \ssp{\int_{0}^{T}\nu(t)g\big(z(t)\big)\,dt+\int_T^\infty\nu(t)\,dt\cdot
   \mathcal{V}\big[\nu^{T}_{\textrm{shift}}
   \big]\big(z(T)\big)}{\omega}-
   \mathcal{V}[\nu](\omega)\bigg|\leq\varkappa.
   \end{equation}

\subsection*{Шаг 2. Выбор и обработка плотности}
Зафиксируем некоторые положительное число
 $\gamma<\frac{\varepsilon}{2M}\min\big(\check{\lambda},\hat{\lambda},1-p,1/2\big)$ и плотность~$\hat{\mu}$,  для которой выполнено
\begin{eqnarray}\label{QQ}
\sup_{t>0}\hat\mu(t)<\gamma,\quad V_{0}^{q[\hat\mu](1-\varepsilon)}[\hat{\mu}]\cdot q[\hat{\mu}](1-\varepsilon)\leq M.
\end{eqnarray}

Зададим отображение $\mu:\mathbb{R}_+\to\mathbb{R}_+$ по следующему правилу:
$$\mu(t)=\hat\mu(t)+\frac{\varepsilon}{q[\hat{\mu}](1-\varepsilon)}\qquad\forall t\in[0,q[\hat{\mu}](1-\varepsilon)]$$
и  $\mu(t)=0$ при $t>q[\hat{\mu}](1-\varepsilon)$.
Легко убедиться, что
\begin{eqnarray*}
\int_0^{\infty}\mu(t)\, dt=\int_0^{q[\hat{\mu}](1-\varepsilon)}\mu(t)\, dt=\int_0^{q[\hat{\mu}](1-\varepsilon)}\hat{\mu}(t)\, dt+\varepsilon=1;
\end{eqnarray*}
таким образом
$\mu$ является плотностью, и  $q[{\mu}](r)<q[\hat{\mu}](1-\varepsilon)$ для всех $r\in(0,1).$

Теперь, в силу
 \begin{equation*}
1-\varepsilon=\int_0^{q[\hat{\mu}](1-\varepsilon)}{\hat\mu}(s)\,ds
\leq\int_0^{q[\hat{\mu}](1-\varepsilon)}\sup_{t\geq 0}\hat{\mu}(t)\,ds=
q[\hat{\mu}](1-\varepsilon)\sup_{t\geq 0}\hat{\mu}(t),
\end{equation*}
 с учетом $M>1,\varepsilon<1/10,$
 мы также имеем
        \begin{eqnarray*}
        \sup_{t\geq 0}\mu(t)\leq
        V_0^\infty [\mu]&\leq&
        V_0^{q[\hat{\mu}](1-\varepsilon)}[\hat{\mu}]+
        \frac{\varepsilon}{q[\hat{\mu}](1-\varepsilon)}\\
        &\leq&
        \frac{M+\varepsilon}{q[\hat{\mu}](1-\varepsilon)}\leq
        \frac{M+\varepsilon}{1-\varepsilon}\sup_{t\geq 0}\hat\mu(t)<2M\sup_{t\geq 0}\hat\mu(t)<2M\gamma,\\
\int_{0}^\infty |\mu(t)-\hat{\mu}(t)|\,dt&=&
\int_0^{q[\hat{\mu}](1-\varepsilon)}\frac{\varepsilon}{q[\hat{\mu}](1-\varepsilon)}\, dt+\int_{q[\hat{\mu}](1-\varepsilon)}^\infty\hat{\mu}(t)\, dt=
2\varepsilon,
\end{eqnarray*}
откуда
    \begin{eqnarray}\label{sup}
\sup_{t\geq 0}\mu(t)\leq V_0^\infty [\mu]<2M\gamma, \quad \int_{0}^\infty|\mu(t)-\hat{\mu}(t)|\,dt\leq 2\varepsilon.
\end{eqnarray}

 Кроме того, отметим, что для всех $x,y>0$ имеет место $$|\ln x -\ln y|=\ln \frac{\max \{x,y\}}{\min \{x,y\}}\leq \frac{\max \{x,y\}}{\min \{x,y\}} -1=\frac{|x-y|}{\min \{x,y\}};$$
теперь из $q[{\mu}](1-\varepsilon)<q[\hat{\mu}](1-\varepsilon)$ и определения $M$, следует, что
$$
V_{0}^{q[{\mu}](1-\varepsilon)}[{\mu}]\leq V_{0}^{q[{\mu}](1-\varepsilon)}[\hat{\mu}]\leq \frac{M}{q[{\hat\mu}](1-\varepsilon)},$$
 откуда
    \begin{eqnarray}\label{log}
     V_0^{q[{\mu}](1-\varepsilon)}[\ln \mu]\leq\frac{V_0^{q[{\mu}](1-\varepsilon)}[\mu]}{\inf_{t\in[0,q[{\mu}](1-\varepsilon))}\mu(t)}\leq
        \frac{M q[\hat{\mu}](1-\varepsilon)}{\varepsilon
q[\hat{\mu}](1-\varepsilon)}=
        \frac{M}{\varepsilon}.
\end{eqnarray}

\subsection*{Шаг 3. Разбиение промежутка}

Разделим полуинтервал $\big[0,q[\mu](1-p^{k^2})\big)$ на~$k^2$ непустых интервалов вида
  $[\tau_{m-1},\tau_{m})$:
для каждого $m=1,\dots,k^2$
 \begin{eqnarray*}
   \tau_0=0,\quad\tau_{m}\rav q[\mu](1-p^m),\\
   \lambda_{m}\rav \frac{\int_{\tau_{m-1}}^{\tau_{m}}
   \mu(t)\,dt}{\tau_{m}-\tau_{m-1}}=\frac{p^{m-1}(1-p)}{\tau_{m}-\tau_{m-1}}\leqref{sup}2M\gamma<\varepsilon\min(\hat\lambda,\check{\lambda},1-p),\\
   \Delta\tau_{m}\rav\tau_{m}-\tau_{m-1}>p^{m-1}/\varepsilon>p^{k^2}/\varepsilon=1.
 \end{eqnarray*}
Отметим, что поскольку для всех
   $m=1,\dots,k^2$ имеет место
 \begin{eqnarray*}
   \int_{0}^{\Delta\tau_{m}}\lambda_m p^{1-m}\,dt=\Delta\tau_{m}\,\lambda_m p^{1-m}=1-p<r_0<\delta_\varepsilon,\\
   {\lambda_m p^{1-m}}<   {2M\gamma p^{-k^2}}=
  2M\gamma/\varepsilon<\min\big(\check{\lambda},\hat{\lambda}\big),
 \end{eqnarray*}
то оценки  \eqref{1157},\eqref{1053} и
  \eqref{105353}  выполнены для  $\nu=\varrho_{\lambda_m p^{1-m}}$  и  $T=\Delta\tau_{m}>1$
   при любых
   $m=1,\dots,k^2$. Теперь,  при любых
   $m=1,\dots,k^2$, в силу $T\leq\frac{\delta_\varepsilon}{\lambda_m p^{1-m}}$
 мы имеем
\begin{eqnarray}\nonumber
   U_*(\omega)&\leqref{1053}&\varkappa+\mathcal{V}[\nu](\omega)\\
   &\leqref{105353}&
      2\varkappa+\ssp{\int_{0}^{\Delta\tau_{m}}\nu(t)g\big(z(t)\big)\,dt+
      \int_{\Delta\tau_{m}}^\infty\nu(t)\,dt\cdot
   \mathcal{V}\big[\nu^{\Delta\tau_{m}}_{\textrm{shift}}
   \big]\big(z(\Delta\tau_{m})\big)}{\omega}\nonumber\\
   &\leqref{1053}&
      3\varkappa+\ssp{\int_{0}^{\Delta\tau_{m}}\nu(t)g\big(z(t)\big)\,dt+
      \int_{\Delta\tau_{m}}^\infty\nu(t)\,dt\cdot
   U_*\big(z(\Delta\tau_{m})\big)}{\omega}\nonumber\\
   &\leqref{1157}&
      3\varkappa+2\varepsilon(1-p)+\ssp{\int_{0}^{\Delta\tau_{m}}\lambda_m p^{1-m}g\big(z(t)\big)\,dt+
      p
   U_*\big(z(\Delta\tau_{m})\big)}{\omega}\nonumber\\
   &=&
      5\varkappa+\ssp{\int_{0}^{\Delta\tau_{m}}\lambda_m p^{1-m}g\big(z(t)\big)\,dt+
      p
   U_*\big(z(\Delta\tau_{m})\big)}{\omega}.\label{1617}
   \end{eqnarray}

\subsection*{Шаг 4. Некорректные промежутки}

Назовем полуинтервал  $[\tau',\tau'')\subset \mathbb{R}$ корректным, если выполнено
 $$V_{\tau'}^{\tau''}[\ln\mu]\leq \frac{M}{k\varepsilon}.$$
В \eqref{log} показано, что $V_0^{q[{\mu}](1-\varepsilon)}[\ln \mu]\leq \frac{M}{\varepsilon}$, тогда   среди всех интервалов $[\tau_{m-1},\tau_{m})$  $(m=1,\dots,k^2)$ не более~$k$ некорректных.

Зададим отображение $\tilde{\mu}:\mathbb{R}_+\to\mathbb{R}$ правилом: $\tilde{\mu}(t)=\mu(t)$
для всех $t\geq \tau_{k^2}$, $\tilde{\mu}(t)=\mu(t)$ для всех $t\in[\tau_{m-1},\tau_{m})$, если этот полуинтервал корректный, и
     $\tilde{\mu}(t)=\lambda_m$ для всех $t\in[\tau_{m-1},\tau_{m})$, если этот полуинтервал некорректный.
Легко проверить, что $\tilde{\mu}\in\mathfrak{D}.$

 Более того, для всех полуинтервалов вида $[\tau_{m-1},\tau_{m})$ интегралы от $\mu$ и $\tilde{\mu}$  равны  $p^{1-m}(1-p)\leq 1-p=\delta.$
 Следовательно, интеграл от модуля их разности равен нулю для каждого корректного промежутка и не превосходит  $2\delta$ для каждого некорректного. Поскольку некорректных промежутков не более $k$,
то интеграл по $\mathbb{R}_+$ от модуля  разности $\hat\mu$ и $\tilde{\mu}$  не превосходит
$2k\delta<2k\cdot{\varepsilon}/{k}=2\varepsilon$ по выбору $k$.
Теперь из  \eqref{sup} следует, что интеграл
по $\mathbb{R}_+$ от модуля  разности $\hat\mu$ и $\tilde{\mu}$ не превосходит $4\varepsilon.$
Наконец, воспользовавшись \eqref{o1}, имеем
    \begin{eqnarray} \label{1214_}
      \mathcal{V}[\hat\mu]\leq\mathcal{V}[\tilde\mu]+4\varepsilon.
    \end{eqnarray}

Отметим, что  $\tilde{\mu}$ было так определено, что для всех  $m=1,\dots,k^2$ выполнено
 $V_{\tau_{m-1}}^{\tau_{m}}[\ln\tilde\mu]<{M}/{k\varepsilon}$,
 теперь, по выбору~$k$, мы последовательно получаем $k\varepsilon\ln(1+\varepsilon)>M$, $1+\varepsilon>e^{{M}/{k\varepsilon}}$ и
\begin{eqnarray*}
    \lambda_m&=&\frac{\int_{\tau_{m-1}}^{\tau_{m}}
    \tilde{\mu}(t)\,dt}{\tau_{m}-\tau_{m-1}}
   \leq
   \sup_{t\in[\tau_{m-1},\tau_{m})}
   \tilde{\mu}(t)\\
   &\leq& e^{V_{\tau_{m-1}}^{\tau_{m}}[\ln\tilde{\mu}]}\inf_{t\in[\tau_{m-1},\tau_{m})}\tilde{\mu}(t)
    \leq (1+\varepsilon)\inf_{t\in[\tau_{m-1},\tau_{m})}\tilde{\mu}(t).
\end{eqnarray*}
 Отсюда, в силу $g\leq 1$, $p^{m-1}-p^m=\lambda_m\Delta\tau_{m}$,  $\varepsilon(1-p)=\varkappa$, для всех $m=1,\dots,k^2$ получаем
\begin{eqnarray}\nonumber
   \int_{0}^{\Delta\tau_{m}} \lambda_{m}g\big(z(t)\big)\,dt&\leq&
    \varepsilon\int_{0}^{\Delta\tau_{m}} \tilde{\mu}(t+\tau_{m-1})\,dt+
   \int_{0}^{\Delta\tau_{m}} \tilde{\mu}(t+\tau_{m-1})g\big(z(t)\big)\,dt\\
    &=&
    \varkappa p^{m-1}+
    \int_{0}^{\Delta\tau_{m}} \tilde{\mu}(t+\tau_{m-1})g\big(z(t)\big)\,dt.
    \label{509}
\end{eqnarray}

Наконец заметим, что $V_0^\infty[\tilde{\mu}]\leq V_0^\infty [{\mu}] $, в силу \eqref{sup}, мы также имеем
   $
    V_0^\infty[\tilde{\mu}]<2M\gamma<\varepsilon\check{\lambda}$.
Отсюда, при любом $m=1,\dots,k^2$ для
$\nu=\tilde\mu^{\tau_{m-1}}_{\textrm{shift}}\in\mathfrak{D}$ имеем
$$V_0^\infty[{\nu}]=\frac{V_{\tau_{m-1}}^\infty[\tilde{\mu}]}{\int_0^{\tau_{m-1}}\tilde\mu(t)\,dt}=p^{1-m}
V_{\tau_{m-1}}^\infty[\tilde{\mu}]<p^{-k^2}V_{0}^\infty[{\mu}]<\frac{2M\gamma}{\varepsilon}<\check{\lambda}.$$
Теперь применяя \eqref{105353} для $T=\Delta \tau_m>1$, $\nu=\tilde\mu^{\tau_{m-1}}_{\textrm{shift}}\in\mathfrak{D}$,
для каждого   $m=1,\dots,k^2$ имеем
   \begin{eqnarray}
   \label{105353__}
    &\ &\!\!\!\! p^{m-1}\varkappa+\ssp{\int_{0}^{\infty}\!\!\tilde\mu(t+\tau_{m-1})g\big(z_{m}(t)\big)\,dt}{\omega}\\
      &=&\!\!\!\!p^{m-1}\varkappa+
   p^{m-1}
      \ssp{\int_{0}^{\infty}\!\!\nu(t)g\big(z_{m}(t)\big)\,dt}{\omega}\nonumber\\
      &\geq&\!\!\!\!\nonumber
      p^{m-1}\ssp{\int_{0}^{\Delta\tau_m}\!\!\!\!\nu(t)g\big(z_{m}(t)\big)\,dt+
      \ssp{\int_{0}^{\infty}\!\!\nu(t+\Delta\tau_m)g\big(z_{m+1}(t)\big)\,dt}{z_m(\Delta\tau_m)}}{\omega}\\
   &=&\!\!\!\!\nonumber
      \ssp{\int_{0}^{\Delta\tau_m}\!\!\!\!\tilde\mu(t+\tau_{m-1})g\big(z_{m}(t)\big)\,dt+
      \ssp{\int_{0}^{\infty}\!\!\tilde\mu(t+\tau_{m})g\big(z_{m+1}(t)\big)\,dt}{z_m(\Delta\tau_m)}}{\omega}.
   \end{eqnarray}

\subsection*{Шаг 5. Прямой ход}
   Зафиксируем произвольное $\omega\in\Omega.$
Для всех $m=1,\dots,k^2$,  в силу
   \eqref{1617} и~\eqref{509} мы получаем
\begin{eqnarray*}
 p^{m-1}{U_*}(\omega) &\leq&5\varkappa  p^{m-1}+
 p^{m-1}\ssp{\int_{0}^{\Delta\tau_{m}}{\lambda_m p^{1-m}}g\big(z(t)\big)\,dt
+p
   U_*\big(z(\Delta\tau_{m})\big)}{\omega}\\
   &=& 5\varkappa  p^{m-1}+\ssp{\int_{0}^{\Delta\tau_{m}}{\lambda_m}g\big(z(t)\big)\,dt
+p^m
   U_*\big(z(\Delta\tau_{m})\big)}{\omega}\\
   &\leq& 6\varkappa  p^{m-1}+\ssp{\int_{0}^{\Delta\tau_{m}}\tilde{\mu}(t+\tau_{m-1})g\big(z(t)\big)\,dt
+p^m
   U_*\big(z(\Delta\tau_{m})\big)}{\omega}.
    \end{eqnarray*}
В частности, для каждого номера  $m=2,\dots,k^2$, для любого процесса $z_{m-1}\in\mathbb{K}$,
имеем
\begin{eqnarray}\label{zmm}
 &\ &\!\!\!p^{m-1}{U_*}(z_{m-1}(\Delta\tau_{m-1}))\\
 &\leq&\!\!\!6\varkappa p^{m-1}+\!
\ssp{\int_{0}^{\Delta\tau_{m}} \!\!\!\tilde{\mu}(t+\tau_{m-1}) g\big(z_m(t)\big)\,dt
+p^m
   U_*\big(z_m(\Delta\tau_{m})\big)}{z_{m-1}(\Delta\tau_{m-1})}.\nonumber
\end{eqnarray}
Подставляя \eqref{zmm} при  $m=2$ в полученное выше для $m=1$, мы имеем
\begin{eqnarray*}
 {U_*}(\omega) &\leq&6\varkappa+
 \sspk{\int_{0}^{\tau_{1}}\tilde{\mu}(t)g\big(z_1(t)\big)\,dt
+p
   U_*\big(z_1(\tau_{1}-\tau_{0})\big)}{\omega}\\
    &\leqref{zmm}&6\varkappa(1+p)+
 \sspll{\int_{0}^{\tau_{1}}\tilde{\mu}(t)g\big(z_1(t)\big)\,dt}\\
&\ &\ +\ssprr{p\ssp{\int_{0}^{\tau_{2}-\tau_{1}}\tilde{\mu}(t+\tau_{1})g\big(z_2(t)\big)\,dt
+p^2
   U_*\big(z_2(\tau_{2}-\tau_{1})\big)}{z_1(\tau_1)}}{\omega}.
    \end{eqnarray*}
  Аналогично, повторяя
 для $m=3,\dots,k^2-1$, и наконец $k^2$, воспользовавшись $\varkappa=\varepsilon(1-p)$, мы получаем
 \begin{eqnarray*}
{U_*}(\omega)
 &\leq&6\varepsilon(1-p)(1+p+\dots+p^{k^2-1})\\
 &\ &
 +\ssplll{\int_{0}^{\tau_{1}}\tilde{\mu}(t)g\big(z_1(t)\big)\,dt+
 \sspll{\int_{0}^{\tau_{2}-\tau_{1}}\tilde{\mu}(t+\tau_{1})g\big(z_2(t)\big)\,dt}}\\
 &\ &
+\dots+\sspl{\int_{0}^{\Delta\tau_{k^2}}\tilde{\mu}(t+\tau_{k^2-1})g\big(z_{k^2}(t)\big)\,dt}
\\
 &\ &\ssprrr{\ssprr{
 \sspr{+p^{k^2}    U_*(z_{k^2}(\Delta\tau_{k^2}))}{z_{k^2-1}(\Delta\tau_{k^2-1})}\dots
 }{z_1(\tau_1)}}{\omega}.
\end{eqnarray*}
 Теперь в силу $U_*\leq 1, p^{k^2}=\varepsilon$, и $g\geq 0$, также показано
\begin{eqnarray*}
{U_*}(\omega)
 &\leq&7\varepsilon+\ssplll{\int_{0}^{\tau_{1}}\tilde{\mu}(t)g\big(z_1(t)\big)\,dt+
 \sspll{\int_{0}^{\tau_{2}-\tau_{1}}\tilde{\mu}(t+\tau_{1})g\big(z_2(t)\big)\,dt}}\\
 &\ &
+\dots+\sspl{\int_{0}^{\Delta\tau_{k^2}}\tilde{\mu}(t+\tau_{k^2-1})g\big(z_{k^2}(t)\big)\,dt}
\\
 &\ &\ssprrr{\ssprr{
  \sspr{+\sspk{\int_{0}^{\infty}{\tilde\mu}(t+\tau_{k^2})
  g\big(z_{k^2}(t)\big)\,dt}{z_{k^2}(\Delta\tau_{k^2})}}{z_{k^2-1}(\Delta\tau_{k^2-1})}\dots
  }{z_1(\tau_1)}}{\omega}.
\end{eqnarray*}
\subsection*{Шаг 6. Обратный ход}
 Применяя \eqref{105353__} для $m=k^2,\dots,1$,
 мы последовательно получаем
\begin{eqnarray*}
{U_*}(\omega)
 &\leq&7\varepsilon+p^{k^2}\varkappa+\ssplll{\int_{0}^{\tau_{1}}\tilde\mu(t)g\big(z_1(t)\big)\,dt+
 \sspll{\int_{0}^{\tau_{2}-\tau_{1}}\tilde\mu(t+\tau_{1})g\big(z_2(t)\big)\,dt}}\nonumber\\
 &\ &
+\dots+\sspl{\int_{0}^{\infty}\tilde\mu(t+\tau_{k^2-1})g\big(z_{k^2}(t)\big)\,dt}
\ssprrr{\ssprr{\sspr{}{z_{k^2-1}(\Delta\tau_{k^2-1})}\dots}{z_1(\tau_1)}}{\omega}\nonumber\\
&=&\!\!\!\!\!\dots=
7\varepsilon+(p^1+\dots+p^{k^2})\varkappa\nonumber\\
&\ &\ \ \ \ \ \ \ \ +\sspll{\int_{0}^{\tau_{1}}\tilde\mu(t)g\big(z_1(t)\big)\,dt+
 \sspl{\int_{0}^{\infty}\tilde\mu(t+\tau_{1})g\big(z_2(t)\big)\,dt}}
\ssprr{\sspr{}{z_1(\tau_1)}}{\omega}\nonumber\\
&\leq&7\varepsilon+(1+p^1+\dots+p^{k^2})\varkappa+\ssp{\int_{0}^{\infty}\tilde\mu(t)g\big(z_1(t)\big)\,dt}{\omega}\nonumber\\
&\leq&8\varepsilon+\mathcal{V}[\tilde\mu](\omega)\\
&\leqref{1214_}&12\varepsilon+\mathcal{V}[\hat\mu](\omega)
\end{eqnarray*}
 для всех  $\omega\in\Omega$ и  для всех плотностей с $\hat\mu$ со свойством \eqref{QQ}.
 
 Предложение доказано.

\end{document}